%% file: marticl1.tex
\documentclass{article} 
\usepackage[all]{xy} \SelectTips{cm}{} 
\input{./mformat.tex}
\def\ADDRESS{\par Institute of Mathematics and Informatics
\par Akademijos 4,
\par LT-08663 Vilnius,
\par Lietuva
\par
\par geva \@ ktl.mii.lt
}%

\begin{document}
%
%
%
\BIBLIOGRAPHY
1! Benabou: Ben>
1! Birkhoff1: Bir>
2! Bourbaki: Bou>
4! BourbakiTop: Bou>
1! Freyd: Fre>
1! Eilenberg: Eil>
2! GahlerGA: Gah>
1! Gahler: Gah>
1! Goguen: Gog>
1! Joyce: Joy>
1! MacLane: Mac>
2! Martin: Mar>
3! Valius: Val>
\ENDED                   
\rf
\centerline{%
\bigfont United sight to an algebraic operations and convergence }
\centerline{by \bf Gintaras VALIUKEVI\v CIUS}
\vskip 12pt 
\noindent {\bf Abstract.}
Algebraic operations and algebraic spaces
was invented to deal with various algebraic problems. Nowadays we can say
that was topologization of algebra. The theory of cathegories became
an abstract theory, and classical algebraic spaces became only an examples
of simplest convergence spaces. In this article we deal with convergences
as instances of multivalued appointments and define the continuity
as property  of commuting squares. This is beginning point for
axiomatisation of multivalued mappings continuity  properties.
Math subj. classification 2000: 46 H05, 54 B30. 

\input{./malgebra.tex}
\input{./mconverg.tex}
\input{./madher.tex}
\input{./mborno.tex}
\vskip 24pt
\REFERENCES
\vskip 24pt
\ADDRESS
\vskip 12pt
2010.04.19
%
\end{document}

%% file: mformat.tex
\font\bigfont=cmr17  
\font\rf=cmr10
\font\bf=cmb10
\newcount\partnumber
\expandafter\newcount\csname number1\endcsname%
\expandafter\newcount\csname number2\endcsname%
\expandafter\newcount\csname number3\endcsname%
\expandafter\newcount\csname number4\endcsname%
\newcount\cnumber
\newcount\counter
\def\Part #1. {%
\global\cnumber=0
\global\advance\partnumber by 1
\vskip 12 pt
\noindent{\bigfont\the\partnumber.
\space
#1}\penalty+10000}
\def\Partnumber#1 {%
\vskip 6pt
\if#1, 
\global\advance\expandafter \csname number\the \cnumber\endcsname by +1
\else 
\if#1: 
\global \advance \cnumber by 1
\expandafter\csname number\the \cnumber\endcsname=1
\else
\if#1.\global\advance \cnumber by -1 
\expandafter\advance\expandafter \csname number\the\cnumber\endcsname by +1
\else
\fi
\fi
\fi
\par\noindent{\bf \the\partnumber%
\counter=0%
\loop\ifnum \counter<\cnumber%
\global\advance \counter by 1 %
.\the\expandafter\csname number\the\counter\endcsname \repeat%
}%
\hskip 5pt}

\long\def\Example. #1\Ended{%
\par
\noindent{\bf Example.\enspace}
{\sl #1}\par\irodymopabaiga
}%
\long\def\Definition. #1\Ended{
\par
\noindent{\bf Definition.\enspace}
{\sl #1}\par\irodymopabaiga
}%

\long\def\Proposition. #1\Proof:{%
\par
\noindent {\bf Proposition.\enspace}
{\sl #1}\par
\noindent{\bf Proof: }
}
\def\WORKING#1! #2: #3\gs{%
\edef\skaitlius{\the\REFER}%
\expandafter\edef\csname #2skaitlius\endcsname{\skaitlius}%
\expandafter\edef\csname\skaitlius tartinis\endcsname{#2}
}%
\newcount\REFER
\def\nonloop#1\repeat{\def\body{#1}\noniterate}%
\def\noniterate{\body\let\next\relax\else\let\next\noniterate\fi\next}%
\def\READS#1>#2\gs{\WORKING#1\gs\def\remainder{#2}}%
\def\BIBLIOGRAPHY#1\ENDED{%
\def\remainder{#1}\REFER=0 %
\nonloop
\advance\REFER by 1
\expandafter\READS\remainder\gs
\ifx\remainder\space\repeat\edef\MAXREFER{\the\REFER}%
}%
\long\def\nuoroda* #1: #2>#3*{
\space [\csname #1skaitlius\endcsname]%
\expandafter\ifx\csname#1visas\endcsname\relax
\expandafter\gdef\csname#1visas\endcsname{#3}%
\fi}%
\def\REFERENCES{\REFER=0
}%
\def\normalref#1{\atskiria#1\gs}
\def\atskiria#1=#2\gs{\def\pagrindinis{#2}{\it #1#2}}
\def\nearrow{\delimiter"3225225}
\def\searrow{\delimiter"3226226}
\def\swarrow{\delimiter"322E22E}
\def\nwarrow{\delimiter"322D22D}
\def\mapright #1 {\mathrel{\smash{\mathop{\longrightarrow}\limits^{#1}}%
\vphantom\longrightarrow}}
\def\maprightun #1 {\mathrel{\smash{\mathop{\longrightarrow}\limits_{#1}}%
\vphantom\longrightarrow}}
\def\mapleft #1 {\mathrel{\smash{\mathop{\longleftarrow}\limits^{#1}}%
\vphantom\longleftarrow}}
\def\mapleftun #1 {\mathrel{\smash{\mathop{\longleftarrow}\limits_{#1}}%
\vphantom\longleftarrow}}

\makeatletter  
\def\setbalta{\setbox\@ne\null
\ht\@ne\ht\z@ 
\dp\@ne\dp\z@ 
\wd\@ne\wd\z@ 
\box\@ne}%
\def\baltalongrightarrow{\setbox\z@\hbox{$\m@th\longrightarrow$}%
\mathrel{\setbalta}}%
\def\baltalongleftarrow{\setbox\z@\hbox{$\m@th\longleftarrow$}%
\mathrel{\setbalta}}%
\def\maprightdouble^#1_#2 {\mathrel{\mathrel{\mathop{\baltalongrightarrow}%
\limits^{#1}_{#2}}%
\mathrel{\llap{\raise 1pt\hbox{$\m@th\longrightarrow$}}}%
\mathrel{\llap{\lower 1pt\hbox{$\m@th\longrightarrow$}}}}}
\def\mapleftdouble^#1_#2 {\mathrel{\mathrel{\mathop{\baltalongrightarrow}%
\limits^{#1}_{#2}}}%
\mathrel{\llap{\raise 1pt \hbox{$\m@th\longleftarrow$}%
\llap{\lower 1pt \hbox{$\m@th\longleftarrow$}}}}}
\def\longrightleftarrow{\mathrel{\mathrel{\raise 1pt
\hbox{$\m@th\longrightarrow$}}\joinrel%
\mathrel{\llap{\lower 1pt\hbox{$\m@th\longleftarrow$}}}}}
\def\longleftrightarrow{\mathrel{\mathrel{\raise 1pt
\hbox{$\m@th\longleftarrow$}}\joinrel%
\mathrel{\llap{\lower 1pt\hbox{$\m@th\longrightarrow$}}}}}
\def\suporavimas^#1_#2 {\mathrel{\mathrel{%
\mathop{\baltalongrightarrow}\limits^{#1}_{#2}}%
\mathrel{\llap{\raise 1pt\hbox{$\m@th\longrightarrow$}}}
\mathrel{\llap{\lower 1pt\hbox{$\m@th\longleftarrow$}}}}}
\def\sansuporavimas^#1_#2 {\mathrel{\mathrel{\mathop{\baltalongrightarrow}%
\limits^{#1}_{#2}}%
\mathrel{\llap{\raise 1pt\hbox{$\m@th\longleftarrow$}}}%
\mathrel{\llap{\lower 1pt\hbox{$\m@th\longrightarrow$}}}}}
\def\rightarrowdouble{\Longrightarrow}
\def\leftarrowdouble{\Longleftarrow}
\def\maprightdble #1 {\mathrel{\smash{\mathop{\rightarrowdouble}\limits^{#1}}%
\vphantom\rightarrowdouble}}
\def\mapleftdble #1 {\mathrel{\smash{\mathop{\leftarrowdouble}\limits^{#1}}%
\vphantom\leftarrowdouble}}
\def\mapup #1 {\Big\uparrow\mathrel{\rlap{$\vcenter{\hbox{$\scriptstyle#1$}}$}}}
\def\mapupl #1 {\llap{$\vcenter{\hbox{$\scriptstyle#1$}}$}\Big\uparrow}
\def\mapdown #1 {\Big\downarrow\rlap{$\vcenter{\hbox{$\scriptstyle#1$}}$}}
\def\mapdownl #1 {\llap{$\vcenter{\hbox{$\scriptstyle#1$}}$}\Big\downarrow }
\def\mapdowndouble #1 {\llap{$\vcenter{\hbox{$\scriptstyle#1$}}$}
\downarrowdouble}%
\def\mapdowndoublek #1 {%
\downarrowdouble\rlap{$\vcenter{\hbox{$\scriptstyle#1$}}$}}%
\def\mapdowndoublelr #1*#2 {\llap{$\vcenter{\hbox{$\scriptstyle#1$}}$}
\downarrowdouble\rlap{$\vcenter{\hbox{$\scriptstyle#2$}}$}}%
\def\mapupdouble #1 {\llap{$\vcenter{\hbox{$\scriptstyle#1$}}$}
\uparrowdouble}
\def\mapupdoublel #1 {%
\uparrowdouble\rlap{$\vcenter{\hbox{$\scriptstyle#1$}}$}}
\def\mapupdoublelr #1*#2 {\llap{$\vcenter{\hbox{$\scriptstyle#1$}}$}
\uparrowdouble\rlap{$\vcenter{\hbox{$\scriptstyle#2$}}$}}
\def\mapupdown #1*#2 {\llap{$\vcenter{\hbox{$\scriptstyle#1$}}$}
\updownarrow\rlap{$\vcenter{\hbox{$\scriptstyle#2$}}$}}
\def\mapdownup #1*#2 {\llap{$\vcenter{\hbox{$\scriptstyle#1$}}$}
\downuparrow\rlap{$\vcenter{\hbox{$\scriptstyle#2$}}$}}
\def\mapne #1 {\Big\nearrow
\rlap{$\vcenter{\hbox{$\scriptstyle#1$}}$}}%
\def\mapnel #1 {\llap{$\vcenter{\hbox{$\scriptstyle#1$}}$}
\Big\nearrow}%
\def\mapneup #1 {\smash{\mathop{\Big\nearrow}\limits^{#1}}}
\def\mapneun #1 {\smash{\mathop{\Big\nearrow}\limits_{#1}}}
\def\mapse #1 {\Big\searrow\rlap{$\vcenter{\hbox{$\scriptstyle#1$}}$}}
\def\mapsel #1 {\llap{$\vcenter{\hbox{$\scriptstyle#1$}}$}
\Big\searrow}%
\def\mapseup #1 {\smash{\mathop{\Big\searrow}\limits^{#1}}}
\def\mapseun #1 {\smash{\mathop{\Big\searrow}\limits_{#1}}}
\def\mapsw #1 {\Big\swarrow\rlap{$\vcenter{\hbox{$\scriptstyle#1$}}$}}
\def\mapswl #1 {\llap{$\vcenter{\hbox{$\scriptstyle#1$}}$}\Big\swarrow}
\def\mapswup #1 {\smash{\mathop{\Big\swarrow}\limits^{#1}}}
\def\mapswun #1 {\smash{\mathop{\Big\swarrow}\limits_{#1}}}
\def\mapnw #1 {\Big\nwarrow\rlap{$\vcenter{\hbox{$\scriptstyle#1$}}$}}
\def\mapnwl #1 {\llap{$\vcenter{\hbox{$\scriptstyle#1$}}$}\Big
\nwarrow}
\def\mapnwup #1 {\smash{\mathop{\Big\nwarrow}\limits^{#1}}}
\def\mapnwun #1 {\smash{\mathop{\Big\nwarrow}\limits^{#1}}}
\def\mgrel{\mathrel>}

\def\injrightarrow{\mgrel\joinrel\longrightarrow}%
\def\mazanearrow{\scriptscriptstyle \nearrow}
\def\mazanwarrow{\scriptscriptstyle \nwarrow}
\def\mazasearrow{\scriptscriptstyle \searrow}
\def\mazaswarrow{\scriptscriptstyle \swarrow}

\def\gimimas{\odot}

\def\duoda{\Longrightarrow}

\mathchardef\vertchar="326A 
\def\ventrightarrow{\mathrel{\rightarrow\joinrel\relbar}}
\def\ventleftarrow{\mathrel{\relbar\joinrel\longleftarrow}}
\def\ventuparrow{\vbox{\offinterlineskip
\halign{##\cr$\mskip 1.5 mu\vertchar$\cr $\uparrow$\cr}\vskip -7 pt}}
\def\ventdownarrow{\hskip 1.5pt\vbox{\offinterlineskip
\halign{##\cr$\downarrow$ \cr$\mskip 1.5 mu \vertchar$\cr}
\vskip -7 pt}\hskip -1.5pt}

\def\ventright #1 {\smash{\mathop{\ventrightarrow}\limits^{#1}}}%
\def\ventrightun #1 {\smash{\mathop{\ventrightarrow}\limits_{#1}}}%
\def\ventleft #1 {\smash{\mathop{\ventleftarrow}\limits^{#1}}}%
\def\ventleftun #1 {\smash{\mathop{\ventleftarrow}\limits_{#1}}}%
\def\ventup #1 {\ventuparrow\rlap{$\vcenter{\vskip 7pt\hbox{$\scriptstyle#1$}}$}}
\def\ventupl #1 {\llap{$\vcenter{\vskip 7pt \hbox{$\scriptstyle#1$}}$}\ventuparrow}
\def\ventdown #1 {\ventdownarrow\rlap{$\vcenter{\vskip 7pt\hbox{$\scriptstyle#1$}}$}}
\def\ventdownl #1 {\llap{$\vcenter{\vskip 7 pt \hbox{$\scriptstyle#1$}}$}
\ventdownarrow }

\def\otline{\mathrel{\lower 1.2pt \hbox{\m@th $\otchar$}}}

\mathchardef\topchar="323E  
\mathchardef\potchar="323F  
\mathchardef\vertchar="326A 

\mathchardef\tustis="001F
\mathchardef\visiems="0238 
\mathchardef\esti="0239

\def\sqr#1#2{{\vcenter{\vbox{\hrule height.#2pt
\hbox{\vrule width.#2pt height#1pt \kern#1pt
\vrule width.#2pt}\hrule height.#2pt}}}}
\def\square{\mathchoice{\sqr66}{\sqr34}{\sqr{2.1}3}{\sqr{1.5}3}}
\def\irodymopabaiga{\vskip -12pt\hfill$\displaystyle\square$\vskip 12pt}%
\def\frac #1/#2 {\leavevmode\kern-.1em \raise.5ex\hbox{$\the\scriptfont0
#1$}\kern-.1em/\kern-.15em\lower.25ex\hbox{$\the\scriptfont0 #2$}}
\def\fract#1:#2 {{#1\over#2}}
\def\Rym#1'{\uppercase\expandafter{\romannumeral#1}}%
\def\duoda{\Longrightarrow}

\def\ts{\advance\displaywidth by-\predisplaysize
\advance\displayindent by\predisplaysize
\abovedisplayshortskip -12pt \belowdisplayshortskip 3pt}
\def\VA{{\cal A}}%
\def\VB{{\cal B}}%
\def\VC{{\cal C}}%
\def\VD{{\cal D}}%
\def\VE{{\cal E}}%
\def\VF{{\cal F}}%
\def\VG{{\cal G}}%
\def\VH{{\cal H}}%
\def\VK{{\cal K}}%
\def\VN{{\cal N}}%
\def\VR{{\cal R}}%
\def\VO{{\cal O}}%
\def\VT{{\cal T}}
\def\VV{{\cal V}}%

\makeatother 

\def\@{{\rm \char'100 }}
\def\identity{{\tt id}}
\def\Identity{{\tt Id}}
\def\inv{^{\mit -1}}
\def\tsk{\ .}
\def\kbl{\ ,\quad}
\def\trik{\bigtriangleup}

\def\isvk#1 {D^{\raise0.25em\hbox{$\scriptstyle#1$}}\kern-.125em }

\def\lim{\mathord{\tt lim}}

\def\Chr{\raise 2pt\hbox{$\displaystyle\Gamma$}\!_ U}
\def\VChr{\raise2pt\hbox{$\displaystyle\Gamma$}!_V}
\def\Chrviens{\raise2pt\hbox{$\displaystyle\Gamma$}\!_1}%
\def\Chrdu{\raise2pt\hbox{$\displaystyle\Gamma$}\!_2}
\def\mathtext#1{\mathchoice{#1}{#1}{\hbox{\scriptsize#1}}{\hbox{\tiny#1}}} 
\def\Ins{\mathord{\tt Ins}}
\def\Int{\mathord{\tt Int}}
\def\Adh{\mathord{\mathtext{\tt Adh}}} 
\def\Sid{\mathord{\tt Sid}}
\def\Set{\mathord{\tt Set}}

\def\virsus#1\apacia#2\viskas{\scriptstyle #1\atop\scriptstyle #2}
\def\Cup{\bigcup}
\def\Cap{\bigcap}

\mathchardef\nuo="3204 
\def\kt#1{\rlap{\tt \char"5C}}
\mathchardef\tustis="001F 
\mathchardef\visiems="0238 
\mathchardef\esti="0239
\def\ldouble{<\kern-1pt<}%
\def\rdouble{>\kern-1pt>}%
\def\liet:#1:{\hbox{\sevenrm #1}}
\def\pavirsiaus{\bigcirc\kern-11pt}
\mathchardef\skyr="313B

\def\s{$\sigma$}

\def\prie#1.{\smash{\mathop{\to}\limits_{#1} }}
\def\skyr{|}
\def\mzl{\leq}

\def\dgl{\geq}

\newcommand{\centre}[1]{\save [].[dr]!C*{(#1)}\restore}
\newcommand{\vent}[1][d]{\ar@{-}[#1]|-{\object@{>}}}


%% file: malgebra.tex
\Part  Algebras and their homomorphisms.
\Partnumber:
\normalref{=Algebra} $\VA$ will be understood as the carrying set $\VA$
with some operations. The \normalref{algebraic =operation} on
the set $\VA$ defined  by Birkhoff
\nuoroda* Birkhoff1: Birk> Birkhoff G. On the structure of abstract algebras,
433--454,
Proc. Cambridge Phil. Soc. {\bf 31} (1935).* defined
is any function
over the  set of $n$-tuples $f:\VA^n\longrightarrow \VA$.
The set of operations names $\Sigma$ is
called the \normalref{=signature of algebra}.

\Example.
In the space $X$ the potential set $\VA=2^X$ is defined
as set of all subsets $A\subset X$. It can be considered as Boole algebra
with such operations:

The least element $0:1\longrightarrow \VA$  is defined
as a function from the singlepoint set $1=\lbrace 0\rbrace$
to the set of all subsets $2^X$, appointing 
the wide set $\tustis\in 2^X$ for the  unique point $0\in 1$.

The largest element $t:1\longrightarrow \VA$ is
defined as a function
appointing the whole
space $X\in 2^X$ for the unique point $0\in 1$.

The complement $c:\VA\longrightarrow \VA$ is defined as a function,
appointing the  complement
$A^c= X\setminus A$ for every partial set $A\subset \VA$.

The intersection is defined as a function
$\land:\VA\times\VA\longrightarrow \VA$, appointing the set
 of common points $A\cap B$
for every pair of partial sets
$\langle A,B\rangle$.

The union is defined as a function $\lor: \VA\times \VA\longrightarrow \VA$
appointing the set of points belonging to either set
$A\cup B$ for every  pair of partial sets $\langle A,B\rangle$.
\Ended

A \s-algebra $\VA$ will be an example of nonclassical algebra. In this algebra
we have countable operations $\VA^\VN\rightarrow \VA$. Usually such
operations are called as convergence of countable sequences.

\Partnumber,
The application $f:X\longrightarrow Y$  will be called \normalref{=morphism} from
the set $X$ to the set $Y$. The first set $X$ will be called a \normalref{=source
space}, and the second set $Y$ will be called a \normalref{=target space} of morphism.
For the point in the source space $x\in X$ the application appoints
the point in the target space $f(x)\in Y$, sometimes we denote
$$y=f(x)\uparrow x\in X\tsk$$

 We have the law of composition for two application with
 \normalref{=intermediant space},
 i.e. the target space of the first application $f:X\longrightarrow Y$
 coincides with the source space of second application $g:Y\longrightarrow Z$.
 The result will be $f\circ g: X\longrightarrow Z$
 an application from the source space of the first application to
 the target space of the second application. For the point $x\in X$ it
 appoints the point $z\in Z$ which is calculated
 $$z=(f\circ g)(x)=g(f(x))\tsk$$  For every set $X$ we
 define the identity application $\Identity_X: X\longrightarrow X$
 which for the point $x\in X$ appoints  the same point.
 It will be neutral for  the law of composition, i.e. the identity
 application doesn't change
 such product, if this product is possible
 $$\Identity_X\circ f=f=f\circ \Identity_Y\tsk$$

Bijective application $f:X\longrightarrow Y$
will have an inverse application $g: Y\longrightarrow X$ for the law of
composition
$$ f\circ g=\identity_X\kbl g\circ f=\Identity_Y \tsk$$
It will be called \normalref{=isomorphism} of sets $X$ and $Y$.

\Partnumber,
We remind the construction of set product.
At first we define the \normalref{=two sets multiplication}.
The multiplication of two sets $X$ and $Y$ is defined
as the set of ordered pairs
$$X\times Y= \lbrace \langle a1,a2\rangle: a1\in X, a2\in Y\rbrace\tsk$$

Two sets multiplication is also defined for morphisms. For any pair
of sets morphisms $f:X\longrightarrow Z$
and $g:Y\longrightarrow W$ we define a morphism between sets multiplications
$$f\times g: X\times Y\longrightarrow Z\times W$$ which
 appoints
$$(f\times g)(\langle x,y\rangle)=\langle f(x),g(y)\rangle\in Z\times W
\uparrow\langle x,y\rangle\in X\times Y\tsk$$
Such multiplication maintains the identity morphisms
$$(\Identity_X\times \Identity_Y)=\Identity_{X\times Y}$$
and the composition of morphisms, i.e. for the pairs of composable
morphisms $\langle f_1,f_2\rangle$ and $\langle g_1, g_2\rangle$
we have the pair of composable
products
$\langle f_1\times g_1, f_2\times g_2\rangle$ and equality
$$(f_1\circ g_1)\circ (f_2\times g_2)=
(f_1\times g_1)\circ (f_2\times g_2)\tsk$$
Therefore the two sets mwltiplication can be concidered as bifunctor.

The sets product can be generated by two sets multiplication.
The $n$-tuple $\langle x1,x2,\dots,xn\rangle\in \VA^n$
can be considered as vector
with $n$
components, i.e. the mapping  from index set $[n]$ to the space of values
$\VA$. Any multiplication  of $n$ sets $A_1$, $A_2$, \dots, $A_n$
can be identified with the set of $n$-tuples with only one  isomorphism
of sets. Therefore any two multiplications  of the same sets
are identified with \normalref{=distinguished isomorphisms} of
sets, i.e. we indicate unique one from possible isomorphisms, which
identifies taken two multiplications. The sets product is defined as
\normalref{=abstract notion} in the category of sets, and sets
multiplication will be only \normalref{=concret presentation}
of such abstract notion.

The singlepoint set $I=\lbrace *\rbrace$ is neutral for
set multiplication, i. e for this set we have distinguished isomorphisms
$l_X:I\times X\longrightarrow X$ and $ r_X: X\times I\longrightarrow X$
which appoint
$$l_X(\langle *,x\rangle)=x\kbl r_X(\langle x,*\rangle)=x\uparrow x\in X\tsk$$
These isomorphisms are \normalref{=natural}
for morphisms of sets $f: X\longrightarrow Y$, i.e. we have the
 commuting diagrams
$$\xymatrix{
     I\times X  \ar[d]_{1\times f}\ar[r]^-{l_X} & X \ar[d]^f \\
          I\times Y \ar[r]_-{l_Y} & Y}
\qquad
   \xymatrix{X\times I \ar[r]^-{r_X} \ar[d]_{ f\times 1} & X \ar[d]^f \\
                     Y\times I \ar[r]_-{r_Y} & Y } $$
Also we have coinciding distinguished isomorphisms for one point set
$I\times I\longrightarrow I$,
$$l_I(\langle *,*\rangle)=*= r_I(\langle *,*\rangle)\tsk$$

Next we indicate the distinguished isomorphism
$a_{XYZ}: (X\times Y )\times Z\longrightarrow X\times (Y\times Z)$
defining
the associative equality
$$a_{XYZ}(\langle\langle x,y\rangle,z\rangle)=
\langle x,\langle y,z\rangle\rangle\uparrow x\in X,y\in Y,z\in Z\tsk$$
 They also are natural for morphisms of sets,
i.e. for morphisms $f:X\longrightarrow  X'$, $g: Y\longrightarrow Y'$,
$h: Z\longrightarrow Z'$ we get commuting diagrams
$$\xymatrix @C=4em{ (X\times Y)\times Z \ar[r]^-{a_{XYZ}^{}} \ar[d]_{(f\times g)}
  & X\times (Y\times Z) \ar[d]^{f\times (g\times h)} 
  \\ (X'\times Y')\times Z' \ar[r]_-{a_{X'Y'Z'}^{}} 
  & X'\times (Y'\times Z')}$$
Such distinguished izomorphisms we got by identification of taken
multiplications with the space of triples $X\times Y\times Z$.
Similar identifications also help to prove the pentagonal
diagram identity
$$\xymatrix @C=-0.8em @R=1.5em{  &  (X\times Y)\times (Z\times W)   
  \\ X\times (Y\times (Z\times W)) \ar[ur]!LD^-a
  & &   ((X\times Y)\times Z)\times W \ar@{<-}[ul]!RD_-a
  \\ \save []+<2.2em,-1em>*+{(X\times (Y\times Z))\times W} ="kairys"
          \ar@{<-}[u]_-{1\times a} \restore
  & & \save []+<-2.2em,-1em>*+{X\times((Y\times Z)\times W)}
              \ar@{<-}"kairys" ^a \ar[u]_-{a\times 1} \restore 
            }$$
In Mac  Lane book 
 \nuoroda* MacLane: Mac> MacLane S. Categories for the working mathematician,
Springer: New York -- Heidelberg -- Berlin 1971.*
the coherence theorems asserts that we can provide the distinguished
isomorphisms
for arbitrary finite sets multiplications defined by different binary trees,
if these trees according with properties
of associative monoid results the same word.
I put so much attention to such usual product of sets,
because it will be an explicit example for another more interesting
products.

If we add the distinguished isomorphisms for twisting
$$t_{XY}: X\times Y\longrightarrow Y\times X \kbl t_{X,Y}(\langle x,y\rangle)
=\langle y,x\rangle
\uparrow x\in X,y\in Y\kbl$$
our product becomes commuting one, and we can get all properties
of commutative monoid. The realization of product with the set of pairs
remains the same, but now we identify two symmetric products with the help
of twisting isomorphisms.

Otherwise we can deal the product of sets without identifying what
singlepoint set
will be neutral. Then for the product of sets we get the properties
of the semigroup, without distinguished isomorphisms $l_X$ and $r_X$.

\Partnumber,
Next we see see another special properties of set product.
Every concret realization of the set product will
be coadjoint functor to diagonal one.

Diagonal functor for every set $X$ apoints the couple of sets $<X,X>$,
and for ervery morphism $f: X\rightarrow Y$ appoints the couple of
morphisms $<f,f>$. For noncommuting product we take these couples
noncommuting
$$<f,g>\not=<g,f>\tsk$$

For every set $X$ we define the diagonal morphism
$$i_X: X\longrightarrow X\times X\kbl i_X(x)=
\langle x,x\rangle\uparrow x\in X\kbl$$
and for pair of sets $\langle Y,Z\rangle$ we define the projections of product
$$p1_{\langle Y,Z\rangle}: Y\times Z\longrightarrow Y
\kbl p1_{\langle Y,Z\rangle}(\langle y,z\rangle)=y\uparrow
\langle y,z\rangle\in Y\times Z \kbl $$
$$p2_{\langle Y,Z\rangle}: Y\times Z\longrightarrow Z\kbl
p2_{\langle Y,Z\rangle }(\langle y,z\rangle)=z \uparrow
\langle y,z\rangle\in Y\times Z \tsk$$
Such morphisms are natural for set morphisms, i.e. for
morphisms $f:X'\longrightarrow X$, $g:Y\longrightarrow Y'$,
$h:Z\longrightarrow Z'$
we have commuting diagrams:
$$\xymatrix{ X \ar[r]^-{i_X} & X\times X
  \\ Y \ar[u]^-f \ar[r]_-{i_Y} & Y\times Y \ar[u]_{f\times f}}
\qquad
 \xymatrix @C=3em{Y\times Z \ar[r]^-{p1_{\langle Y,Z\rangle}}
 \ar[d]_{g\times h}
   & Y \ar[d]^g
   \\ Y'\times Z'\ar[r]_-{p1_{\langle Y',Z'\rangle}} & Y'}
 \qquad
 \xymatrix @C=3em{Y\times Z  \ar[r]^-{p2_{\langle Y,Z\rangle }}
  \ar[d]_{g\times h}   & Z \ar[d]^h
   \\ Y'\times Z' \ar[r]_-{p2_{\langle Y',Z'\rangle}} & Z'} $$
The first natural morphism $i_X$ defines  mapping over set of
all morphism pairs
$$\phi_{X\langle Y,Z\rangle}:\langle (X,Y),(X,Z)\rangle\longrightarrow
(X,Y\times Z)\kbl$$
which for the pair of morphism $f: X\longrightarrow Y$,
$g: X\longrightarrow Z $ appoints composition of morphisms
$i_X\circ f\times g:X\longrightarrow Y\times Z$. The second
natural morphism pair
$\langle p1_{\langle Y,Z\rangle},p2_{\langle Y,Z\rangle}\rangle$
defines mapping over the set of morphisms
$$\psi_{X\langle Y,Z\rangle}:(X,Y\times Z)\longrightarrow \langle
(X,Y),(X,Z)\rangle\tsk$$
These mappings are natural for the source and target turnings of morphisms,
i.e. for earlier morphisms we have commuting diagrams of mappings
$$ 
\xymatrix @C=3.8em{\langle \langle X,Y\rangle,\langle X,Z\rangle\rangle
\ar[r]^-{\phi_{X\langle Y,Z\rangle}}
\ar[d]_{\langle \langle f,g\rangle , \langle f,h\rangle \rangle}
  & \langle X,Y\times Z\rangle  \ar[d]^{\langle f,g\times h\rangle }
  \\ \langle \langle X',Y'\rangle,\langle X',Z'\rangle \rangle
  \ar[r]_-{\phi_{X'\langle Y',Z'\rangle}}
  & \langle X',Y'\times Y'\rangle }
\qquad
\xymatrix @C=3.8em{\langle X,Y\times Z\rangle
\ar[r]^{\psi_{X\langle Y,Z\rangle}}
  \ar[d]_{(f,g\times h)}
  & \langle \langle X,Y\rangle ,\langle X,Z \rangle \rangle
  \ar[d]^{\langle \langle f,g\rangle,\langle f,h\rangle\rangle}
  \\ \langle X',Y'\times Y' \rangle  \ar[r]_-{\psi_{X'\langle Y',Z'\rangle}}
  & \langle\langle X',Y'\rangle, \langle X',Z' \rangle \rangle }
$$
The natural mappings $\phi_{X,\langle YZ\rangle}$ will define
 \normalref{=adjunction},
and $\psi_{X,\langle YZ\rangle}$ will define
\normalref{=coadjunction}.

The pair of identities for natural morphisms $i_X$ and
$\langle p1_{\langle Y,Z\rangle },p2_{\langle Y,Z\rangle }\rangle$
$$i_X\circ p1_{\langle X,X\rangle }=\Identity_X \kbl
i_X\circ p2_{\langle X,X\rangle }=\Identity_X $$
provides
the \normalref{first identity of =duality}
$$\phi_{X\langle Y,Z\rangle}\circ \psi_{X \langle Y,Z\rangle}=
\Identity_{\langle\langle X,Y\rangle ,\langle X,Z\rangle \rangle}\kbl $$
and the identity
$$i_{Y\times Z}\circ p1_{YZ}\times p2_{\langle Y,Z\rangle}=
 \Identity_{Y\times Z}$$
 provides the \normalref {second identity of =duality}
$$\psi_{X\langle Y,Z\rangle}\circ \phi_{X\langle Y,Z\rangle}=
\Identity_{\langle X,Y
\times Z\rangle }\tsk$$

In the case of such identities  adjunction and coadjunction are reciprocal
isomorphic mappings. Any product realization with such property is
isomorphic with identification
of unique isomorphism. This can be consequence of the theory for
the \normalref{=dual
functors}. We conclude that the multiplication of sets
is presentation
of unique abstract product. This abstract product is called
\normalref{=Cartesian product} of the sets.
Every concrete product will be
an \normalref{=implementation} of such abstract product.
The abstract product is identified with possible unique
change of implementation. The commuting abstract product of the sets
has more changes of implementation, therefore it is more abstract than
noncommuting one.
Freyd and Scedrov \nuoroda* Freyd: Fre> Freyd P. Scedrov A. Categories and allegories, North
Holland: Amsterdam -- New York -- Oxford -- Tokyo 1990.*
use the commuting Cartesian product, Mac Lane in
 \nuoroda* MacLane: Mac> MacLane S. Categories for the working mathematician,
Springer: New York -- Heidelberg -- Berlin 1971.*  has'nt understood
what abstract categorical notions mean, and he makes no explicit
differences between various notions of Cartesian products.
Commuting abstract product can be comfortably defined as limit of discrete
diagram.

For abstract product we have much more
isomorphic implementations. For example we can identify arbitrary
singlepoint set with $I$ using the  unique isomorphism
$f:\lbrace x\rbrace \longrightarrow \lbrace\ast\rbrace$ which appoints
$f(x)=\ast$, therefore such singlepoint set $\lbrace x\rbrace $
becomes neutral for such abstract set
product implementation.

The abstract commuting Cartesian product extends the
operator of intersection in
Lattice theory and for it
the same sign  $X\land Y$ can be applyed.
\Partnumber,
For more complex algebraic operations many sorted algebra can be defined
\nuoroda* Goguen: Gog> Goguen J. Thatcher J. Wagner E. An initial algebra
approach to the specification, correctness, and implementation
of abstract data types, 80--149, Current trends in programming methodology
vol. 4 Data structuring, ed. Yeh R., Prentice Hall: Englewood Cliffs,
New Jersey in USA 1978.*. Let the carrying sets are indexed by the points
from sort set $s\in S$ and the name of operation indicates what sort
of carrying spaces must be taken, i.e. the signature $\Sigma$
is indexed over the space $S^*\times S$, where $S^*$ denotes the
set of all finite words of alphabet $S$. Therefore every operation will
be a function
$$f_{s1,s2,\dots, sn;s}: \VA_{s1}\times \VA_{s2}\times\dots\times \VA_{sn}
\longrightarrow \VA_s\tsk$$

\Example.
A directory $\VC:\VO\longrightarrow \VO'$ is understood as a set of
arrows from some point in the first set $\VO$ to a point in the second
set $\VO'$.
Points in the sets $\VO$, $\VO'$  are understood as vertexes of arrows.
The arrows between two vertexes $X\in\VO$, $Y\in\VO'$
compound an arrows set $\VC(X;Y)$. The first vertex $X$
is called a source
of arrow, the second vertex $Y$ is called a target of arrow.
A source  is an appointment from the set of arrows to the set
of vertexes $s:\VC\longrightarrow \VO$, for each arrow  it appoints
the source vertex of taken arrow.
The target is an appointment from  the
set of arrows to the set of vertexes $t:\VC\longrightarrow \VO'$,
for each arrow it appoints the  target vertex of taken arrow.

The directory is a sorted algebra with sort set $\bf 3$
$$\VA_1:=\VA\kbl \VA_2:= \VO \kbl\VA_3:= \VO\kbl $$
and two operations $f_{1;2}=s:\VA\longrightarrow \VO$ and
$f_{1;3}=t:\VA\longrightarrow \VO$.

The category $\VC$ will be a dirrectory with the same source and target
set $\VO$. Additionally we must choose the unit arrows $1\in (x,x)$
and define the arrows composition
$$f\circ g\in (x,z)\uparrow f\in (x,y), g\in (y,z)\tsk$$
\Ended

\Partnumber,
For more general algebras on the carrying set $A$
we need to chose other interesting products $A^n$,
and then to define algebraic operations $A^n\longrightarrow A$
as morphisms over such product.
\Example.
Let we have two instances of directories with some intermediant vertex
space $\VD:\VO_1\longrightarrow \VO_2$ and $\VD':\VO_2\longrightarrow \VO_3$.
The product of two directories will be a product bundle for the target
and source  appointments
$$\VD\otimes_{\VO_2} \VD':\VO_1\longrightarrow \VO_3\tsk$$
It will be a part of product of arrow sets and can be
defined as equalizer for target and source appointments,
i.e. it will be the set of arrow pairs for which
two maps coincide
$$\VD\times_{\VO_2} \VD'=(t=s')=\lbrace \langle f,g\rangle\in \VD\times \VD'
:t(f)=s'(g)\rbrace\tsk$$
This new \normalref{=directory product} get distinguished isomorphisms
from semilattice of set
product. This isomorphisms are commuting with identity morphisms of
factor spaces, therefore the equalize of target and source appointments
is maintained, and we can take the trace of distinguished isomorphisms over
the set of equalizer. There aren't projections or diagonal morphisms
for such directory product.

The singlepoint set can't be taken as neutral element for directory product.
The neutral element will be a set of vertexes.
We have the distinguished isomorphisms for equalities
$$l_\VD:\VO_1\otimes_{\VO_1} \VD\longrightarrow \VD\kbl
r_\VD: \VD\otimes_{\VO_2} \VO_2\longrightarrow \VD\tsk$$

For the sake of clarity we express the product bundle as the sum of
nonintersecting arrow sets
$$\VD= \Cup\lbrace (a,b): a\in \VO_1, b\in\VO_2\rbrace \kbl
\VD'=\Cup\lbrace (a',b'): a'\in \VO_2, b'\in \VO_3\rbrace\tsk$$
The product of two directories is get by the product of slices
$$\VD\otimes_{\VO_2}\VD'=\Cup\lbrace (a,b)\times (b,c): a\in \VO_1,
c\in \VO_2 \rbrace\tsk$$
The set of vertexes $\VO$ is understood as directory with identity morphism
$\Identity:\VO\longrightarrow \VO$ taken as source and target appointments.
There is exactly one arrow over each diagonal pair of vertexes
$\langle c,c\rangle$.
The products by such directory are defined
$$\VO_1\otimes_{\VO_1} \VD=
\Cup\lbrace \ \lbrace a\rbrace\times (a,b): a\in \VO_1, b\in\VO_2\rbrace \kbl
$$$$\VD\otimes_{\VO_2} \VO_2 =\Cup \lbrace (a,b)\times \lbrace b\rbrace:
a\in\VO_1, b\in\VO_2\rbrace\tsk $$

Therefore the distinguished isomorphisms is defined over slices identifying
the  set product by singlepoint set with the taken set of arrows
$$\lbrace a\rbrace\times (a,b) \longrightarrow (a,b)\kbl
(a,b)\times \lbrace b\rbrace \longrightarrow (a,b)\tsk$$
\Ended

\Partnumber:
The set product provides a new algebra with product operations.
\Proposition.
If the set product interchange with the product of algebra
by isomorphisms
$$X^n\times Y^n \longrightarrow (X\times Y)^n\kbl$$
then the set product of carrying
spaces provides a new algebra with the set product operations.
\Proof:
It is enough to use the inverse of distinguished isomorphisms
$$I\times I\longrightarrow I \kbl
X^n\times Y^n\longrightarrow (X\times Y)^n$$
to define the $0$-degree and $n$-degree operations in $X\times X$
$$\xymatrix{ I\times I \ar[d]_{f_1\times f_2} & I \ar[l] \ar[d]
  \\ X\times Y   \ar[r]_-{1} & X\times Y}
 \qquad
 \xymatrix{ X^n\times Y^n \ar[d]_{f_1\times f_2}
   & (X\times Y)^n \ar[l] \ar[d]^{f}
   \\ X\times Y \ar[r]_-{1} & X\times Y } $$
\irodymopabaiga
Such algebra will be called a \normalref{=set product algebra}.
In the case when the set product is associative and commuting one,
and algebra's operations can be defined over the set product.

For any pair of directories $\VD_1$ and $\VD_2$ the set product
$\VD_1\times \VD_2$ provides a new product directory with product operations, i.e.
the source and target appointments is defined by products
$$s_1\times s_2:\VD_1\times\VD_2\longrightarrow  \VO_1\times\VO_2\kbl
t_1\times t_2: \VD_1\times \VD_2\longrightarrow \VO_1'\times \VO_2'\tsk$$

\Example.
For the categories the set product interchanges
with bundle product, used to define
the operations of category,
$$(\VA\otimes_{\VO}\VA)\times (\VB\otimes_{\VO'} \VB)\longrightarrow
(\VA\times \VB)\otimes_{\VO\times \VO'}(\VA\times \VB)\kbl$$
therefore we can define the appointment of unit arrow
$\VO\times \VO'\longrightarrow \VA\times \VB$ and appointment of arrows
composition $(\VA\times\VB)\otimes_{\VO\times\VO'}(\VA\times\VB)
\longrightarrow \VA\times\VB$. The set product of categories
provides product category with product operations.
\Ended

\Partnumber.
For the algebras $X$ and $Y$ with $1$-degree operation it is easy to define
the sum of such algebras $X+Y$. We take the nonintersecting union of sets
$X+Y$ and define the sticked operation as the sum
 $$f+f':X+Y\longrightarrow X+Y\tsk$$
For intersecting carriers $X\cap Y\not=\tustis$ we must check the
coincidence of sticked operations over the commune part $X\cap Y$.

For $0$-degree operations we need another construction of new operation
over the sum of carriers. The same problem will also arise for
bigger $d\ge 1$ degree operations. For the sum of groups $X+Y$
we must demand
that the units of both groups coincided and decide what would  be
the composition of members
from different sets. The unique solution is to use the product group
compounded of pairs $X\times Y=\lbrace \langle x,y\rangle: x\in X,
y\in Y\rangle\rbrace $ with
the inclusion of summands as partial groups $X\longrightarrow X\times Y$,
$Y\longrightarrow X\times Y$ taking pairs with unit of another sumamnd
$\langle x,e\rangle\in X\times Y$ and $\langle e,x\rangle\in Y$.

The sum of directories poses no problem. The sum of two directories
$\VD_1:\VO_1\longrightarrow \VO_1'$ and $\VD_2':\VO_2\longrightarrow \VO_2'$
will have the sum of arrow sets $\VD_1+\VD_2$ and
will have sourse set the sum $\VO_1+\VO_2$ and target set the sum
$\VO_2 +\VO_2'$. Therefore we can define source appointment
$\VD_1+\VD_2\longrightarrow \VO_1+\VO_2$ and target appointment
$\VD_1+\VD_2\longrightarrow \VO_1'+\VO_2'$.

The sum of categories $\VA+\VB$ will have the same arrows as the sum of
directories. The unit arrow can be taken in each subcategory  separately.
The composition is defined only for arrows from one subcategory in order
to have
common middle vertex.  We shall say that taken categories are
\normalref{=nonintermittent subcategories}
$\VA\subset \VA+\VB$ and $\VB\subset \VA+\VB$.

\Partnumber,
The \normalref{=semigroup} is an example of algebra with one law of multiplication
$$(\circ): X\times X\longrightarrow X\tsk$$
Such algebra is called \normalref{=monoid without unit}.

We ask that this multiplication would be associative, i.e. for all
points we have identities
$$(a\circ b)\circ c= a\circ (b\circ c)\tsk$$

This algebra is not free, we shall say that this algebra
satisfies relations.

The \normalref{=monoid} additionally demands to take the
neutral element $e\in X$ with identities
$$x\circ e=x\kbl e\circ x=x\uparrow x\in X\tsk$$
Such element is understood as operation of $0$-degree
$$e:X^0=I\longrightarrow X\tsk$$

The \normalref{=group} demands to have the inverse element $x\inv\in X$
for every point $x\in X$ with two identities
$$x\circ x\inv=e=x\inv\circ x\uparrow x\in X\tsk$$
The inverse element can be understood as operation of $1$-degree
$$X\longrightarrow X\kbl x\inv\uparrow x\in X\tsk$$
In the case of existence, such element is unique in associative monoid.
For two points $y$ and $z$ wich are inverse elements of $x$ we get
equality
$$y=y\circ e= y\circ (x \circ z)=(y\circ x)\circ z=e\circ z=z\tsk$$
Therefore the inverse element $x\inv$ can be defined as solution of equation
in associative monoid
$$-\circ  x =e\kbl x\circ - =e\tsk$$

\Partnumber,
The category $\VC$ will be a directory with the same source and target
vertexes set $\VO$ and being associative monoid for the directory product.
The multiplication of arrows is called composition
$$(\circ): \VC\otimes_{\VO}\VC \longrightarrow \VC \kbl$$
and neutral element will be appointment of unit arrows
$$1: \VO\longrightarrow \VC\tsk$$

Group for such directory product has its own name
\normalref{=grupoid}. The semigroup for the directory product is not so usual.

\Example.
The class of all directories provides an
example of category with taken composition
and unit arrows.  If we take only directories with finite number of
arrows, we shall not have the units in such category.
Therefore it can be called the category without units.
This provides a serious example of semigroup for the
directory product.
\Ended

We can define the category as an algebra  with operations being
partial mappings over the set  product. Such definition is given in
Freyd and Scedrov
\nuoroda* Freyd: Fre> Freyd P. Scedrov A. Categories and allegories, North
Holland: Amsterdam -- New York -- Oxford -- Tokyo 1990.*
 p. 3. Such algebras with partial
operations they have called  essentially algebraic,
in the sense that we can
construct special products and define functional operations over these
product to get the same (isomorphic) algebra.
In general case the partial operations must be deal as arbitrary
convergence.

\Partnumber,
The morphism between two algebras with the same signature
for set product is defined as an application
$u:\VA\longrightarrow \VB$ which maintains the homologic operations,
 i.e. for each operation name we have the commuting diagram:
$$\xymatrix{\VA^n \ar[r]^-{u^n} \ar[d]_{f} & \VB^n \ar[d]^{f'} \\
          \VA \ar[r]_-{u} & \VB } $$
Others abstract products are defined in its own categories, and
application is changed by arrows between carrying spaces.

\Example.
The functor will be morphism of categories $F: \VA\Longrightarrow \VB$.
At first it must be morphism of underlying directory, i.e. we
must have the application of arrows $F:\VA\longrightarrow \VB$
and application of vertexes $F_0:\VA_0\longrightarrow \VB_0$
which must commute with the source and target appointments. Secondly
we must demand the maintenance of arrow operations over the directory
products. It must maintain the units and the composition of arrows
$$\xymatrix{\VO \ar[r]^-{F_0} \ar[d]_{1} &\VO'\! \ar[d]^{1} \\
          \VA  \ar[r]_-{F} &\VB }
\qquad
\xymatrix @C=4em{\VA\otimes_\VO\VA \ar[r]^-{F\otimes_\VO F} \ar[d]_{(\circ)}
  & \VB\otimes_{\VO'}\VB \ar[d]^{(\circ)}
  \\ \VA \ar[r]_-{F} & \VB } $$
\Ended

All set mappings $f:X\longrightarrow Y$ compounds the
\normalref{=category of sets}, we shall note it
$\Set$.
Vertexes will be all possible (small) sets. The unit arrow will be
identity application $\Identity_X: X\longrightarrow X$. The arrow composition
is taken the usual composition of functions $f\circ g$.

The product category $\Set\times\Set$ has pairs of sets
$\langle X,Y\rangle$ as vertexes
and pairs of mappings $\langle f:X\longrightarrow Z, g:Y\longrightarrow W\rangle$
as arrows. The bifunctor will be a functor over the product category.
The sets multiplication
is an example of bifunctor in the category of sets $\Set$
$$(\times): \Set \times\Set\Longrightarrow \Set\tsk$$

\Partnumber,
More generally operation in the set category
may be arbitrary \normalref{=multivalued function}
$F:\VA^n\ventrightarrow \VA$. Such algebras  are called \normalref{=multialgebras}.
When the operations are arbitrary partial functions, such algebras
are called \normalref{=partial algebras}. If we take
the operations over infinite products $\VA^\VN$, such algebras will be
better interpreted as convergence space. It is nothing with any finite
algebraic
property, but traditionally it will be called algebra. Such tradition is
maintained also by categorical framework of algebras in Mac Lane's book
\nuoroda* MacLane: Mac> MacLane S. Categories for the working mathematician,
Springer: New York -- Heidelberg -- Berlin 1971.* part VI. The united
sight to an algebras calls all algebraic or nonalgebraic operations
as convergences and algebras identifies with topological spaces in
generalized sense as stuctures defined by Bourbaki \nuoroda* Bourbaki: Bou>
Bourbaki N. Th\'eorie des ensembles, Hermann: Paris 1963, Mir: Moskva 1965.*

The category was introduced as ``abstract nonsense'' by S. Eilenberg and
S. MacLane in \nuoroda* Eilenberg: Eil> Eilenberg S. MacLane S.
General theory of natural equivalencees, 231-294, Trans. AMS {\bf 58} (1945).
*. Their theory emphasizes the work with arrows. On the points of
sets constructed ``concrete''
functions are changed by arrows and many properties of arrows can be proved
without investigation of concrete structure of these functions. We can
say that a category is one step of abstraction. We use various structure
in the category's set of arrows, and the second step of abstraction would
be changing this concrete structure in the set of arrows
by some other categorical properties. In such way we encounter the
\normalref{=bicategories}. We want emphasize that categorical
argument are not
absolute fashion of modern mathematics. The ``concrete'' structures explored
by N. Bourbaki must also be admired by ``working mathematicians''.

The algebras is an example of concrete structure. In contrast
the abstract notion of Cartesian product will be an example of
abstract notion. This contrast is very clear in computer programming.
Algebras will be identifyed  with working programs in concrete computer,
and categorical abstract notions will be only specification of working
programs. Nevertheless the categorical approach made essential
progress in todays
programming business.

The monoidal category is defined as category with some bifunctor
$$\square:\VC\times\VC\longrightarrow \VC\tsk$$
It was first explicitly mentioned in 1963 by B\'enabou
\nuoroda* Benabou: Ben> B\'enabou J. Categories avec multiplication,
1887--1890, Comptes Rendue Acad. Sci. Paris {\bf 256} (1963).*.
The name "monoidal" is due to Eilenberg.

This bifunctor is called a multiplication and will
generalize the set multiplication of the set category.
One demands that this multiplication would be associative
after identification
with a natural transformation of functors
$$\alpha: (-\square-)\ \square-\longrightarrow -\square\ (-\square -)\tsk$$
For the coherence one needs that the pentagon diagram
with wedges get from natural transformation $\alpha$ would be commuting.

At this time in physics is introduced the premonoidal categories,
see W. Joyce
\nuoroda* Joyce: Joy> Joyse W. Braided premonoidal Mac Lane coherence,
155--176, J. of Pure and Applied Algebra {\bf 190} (2004).*.
%
%
It is interesting that in his work noncommuting pentagon diagram
are completed to commuting one with a deformation arrow.
The calculating of this arrows provides
nonassociative statistics in quantum physics.

\Partnumber,
A multivalued mapping between points of two sets is identified
with partial set of set multiplication $R\subset X\times Y$. Again the first
space is called a source space
and the second is called a target space. We shall also say that we have
a \normalref{=reform} from the source space to the target space.

The pair of points $\langle x,y\rangle\in R$ is called \normalref{=related}.
More exactly we shall say
that the first point $x\in X$ reforms himself to the second point
$y\in Y$ and the
second point is reformed from the first one.
We can also say that second point is \normalref{=related}
to the first one, and the first
point is \normalref{=corelated} to the second one.
It is more usual to say that the second point $y\in Y$
 is related to the first point $x\in X$ by the relation $R$, or the
 first point $x\in X$ has related point in the target space. But relation
has no source or target space, therefore such language can't be correct
for educated mathematician.

The reform
$A\subset \VO\times\VO'$ will be
an example of directory which has no more than one arrow
for every pair of vertexes.

The reform is uniquely defined by its \normalref{=direct appointment}
$$x\longmapsto x\circ R\kbl X\longrightarrow 2^Y\kbl $$
or by its \normalref{=opposite appointment}
$$y\longmapsto R\circ y\kbl Y\longrightarrow 2^X\tsk$$
The first appointment to a point $x\in X$ appoints all related points
in the target space $x\circ R$. The second appointment to a point $y\in Y$
appoints all the points the source space $x\in X$ to which the taken point
is related.

Two reforms  $F:X\longrightarrow Y$ and $G:Y\longrightarrow Z$
can be composed
$$F\circ G:=\lbrace \langle x,z\rangle\in X\times Z: (\exists y\in Y)
(\langle x,y\rangle\in F,
\langle y,z\rangle\in G)\rbrace\tsk$$

The identity reform is defined by diagonal in set product
$$\Identity_X=\trik\subset X\times X\tsk$$
We can check the associativity for composition of reforms,
and all reforms between arbitrary small sets compound a new
category which will be called the allegory of sets.
It extends the category of sets compounded by all mappings
between arbitrary small sets. 

The reform $F\subset X\times Y$ is called \normalref{=entire} if every point
in the source space has related point in target space.

The reform $F\subset X\times Y$ is called \normalref{=accurate}
if every point in the source space has no more than one related
point in the target space.

Also we rename these properties of opposite reform.

The reform $F\subset X\times Y$ is called \normalref{=covering} if every
point in the target space is related with some point in the source space.

The reform $F\subset X\times Y$ is called \normalref{=exact} if
every point in the target space is related with no more than one point
in the source space.

The entire accurate  reform is a function, the covering function is
\normalref{=surjection},
the exact function is \normalref{=injection}.

For a reform $F:X \longrightarrow Y$ between two spaces
we construct a function of \normalref{=direct image} between potential spaces
$F^P:2^X\longrightarrow 2^Y$
defining \normalref{=image} of arbitrary set in a source space $A\subset X$
$$F^P(A):=\lbrace y\in Y: (\exists x\in A)(\langle x,y\rangle\in F)
\rbrace\kbl$$
and a function of \normalref{=opposite image}
$F^{*P}:2^Y\longrightarrow 2^X$ defining \normalref{=opposite image} of
arbitrary set in a target space $B\subset Y$
$$F^{*P}(B):=\lbrace x\in X: (\exists y\in B)(\langle x,y\rangle\in F)
\rbrace\tsk$$
For the mappings $f:X\longrightarrow Y$ the opposite image is called
an inverse image and is denoted more simply
$$f\inv(B)=f^{*P}(B)\tsk$$

\Partnumber,
The reforms can be abstractly defined in arbitrary regular category, cl.
Freyd and Scedrov
\nuoroda* Freyd: Fre> Freyd P. Scedrov A. Categories and allegories, North
Holland: Amsterdam -- New York -- Oxford -- Tokyo 1990.*.
Cartesian category has the \normalref{=beginning} of every finite diagrams.
It is enough
to have \normalref{Cartesian =product}, \normalref{=terminal object} and
\normalref{=equalizer} of every morphism pair.

All these constructions are \normalref{=abstract}, every  object is
identified only with \normalref{=uniquely existent isomorphism}.
 Regular category additionally demands the existence
of \normalref{=image} for every arrow and the property of product bundle
to maintain
covers, i.e. for covering base mapping we get covering induced mapping
of bunde.  The reform $R: X\ventrightarrow Y$ in Cartesian
category is defined as \normalref{=monopair} of
\normalref{=tabulating mappings}
 $f:Q\longrightarrow X$ and $g: Q\longrightarrow Y$.
The relation is understood as part of Cartesian product defined by this
monopair. It can be named as \normalref{=graphic} of reform
$$G(R)\injrightarrow X\land Y\tsk$$
The reform $R$ can be restore from it's graphic $G(R)\injrightarrow X\land Y$
by choosing what the space is a source and what the space is a target, i.e.
the reform can be denoted by the triple
$$ R=\langle X,G(R),Y\rangle\tsk$$
In regular category we have possibility to  define the composition of two
reforms with an intermediate  space.

The mapping $f:X\longrightarrow Y$
are identifyed with
tabulation $1:X\longrightarrow X$, $f:X\longrightarrow Y$ and arbitrary
reform is equal to fraction of its tabulation $R=f\inv \circ g$.
All reforms of regular category $\VC\subset \VR(\VC)$
compound an allegory. We shall call it a fractional extention.
Functor of Cartesian category will be called representation of
Cartesian category if it maintains
the beginnings of finite diagrams. Functor of regular category will be
called representation of regular category if it additionally maintains
the covers. Every representation of regular category
$F:\VC\Longrightarrow \VC'$ has fractional extension which becomes a
representation of \normalref{=unitary allegories}
$\VR(F):\VR(\VC)\Longrightarrow \VR(\VC')$.

In Cartesian category we can freely define the composition of function
with arbitrary relation $f\circ R$. A topos will be Cartesian category $\VC$
having
potential objects. For any object $C$ we define potential $[C]$
with distinguished reform $\ni: [C]\ventrightarrow C$ wich must  be
initial among reforms to this object by transformations defined
with arrows of category, i.e. for every reform $R: X\ventrightarrow C$
we can find a unique representing arrow $f: X\longrightarrow C$ which
transforms the distinguished reform to the taken one
$$f\circ \ni=R\tsk$$

The set category is an example of topos.
The potential of the set $C$ will be the set of all
subsets $[C]:=2^C$ and distinguished reform every partial set $A\subset C$
relates with its points $x\in C$, we shall write $A\ni x$.

For the reform  $R:X\ventrightarrow C$ the unique representing mapping
$f: X\longrightarrow [C]$ is a direct appointment of this relation.
For any reform $R: X\ventrightarrow Y$ in topos we can define the unique
mapping $R^P:[X]\longrightarrow [Y]$. In the Set category this mapping
coincides with
direct image appointment.
The direct image appointment defines endofunctor of set category
$F: \Set \Longrightarrow \Set$. Its fractional extention is
representation of unitary allegories
$\VR(F):\VR(\Set)\Longrightarrow \VR(\Set)$. We remark that it differs
from earlier direct image appointment.
For the relation $R: X\ventrightarrow Y$ tabulated with the pair of mappings
 $\langle f,g\rangle: Q\longrightarrow
 \langle X, Y\rangle$ the
 set $A\subset X$ has related set $B\subset Y$ exactly  when
there is a set $C\subset Q$ with needed direct images
for tabulating mappings
$$f^P(C)=A\kbl g^P(C)=B\tsk$$
We get maintenance for composition $F\circ G$ and inverce $F^*$
of relations
$$\VR(F\circ G)=\VR(F)\circ \VR(G)\kbl \VR(F^*)=\VR(F)^*\tsk$$
The direct image $f^P(A)\subset Y$ are related to taken set $A\subset X$
by the trace of set product $(A\times f^P(A))\cap R\subset X\times Y$.
It is the biggest of related sets to taken one.


%% file: mconverg.tex
\Part Reform's continuity in convergence spaces.
\Partnumber:
The operations in multialgebras better to interpret as convergences.
Then  morphisms of multialgebras will be defined as mappings
with various continuity properties.

The \normalref{=convergence} in the \normalref{=world space} $X$
is understood as a reform
$\rho\subset X^\Phi\times X$ from sequences
$\VF\in X^\Phi$ to the points
 $x\in X$ of taken set $X$, see G\"ahler
\nuoroda* Gahler: Gah> G\"ahler W. A unified general theory of
convergence, 81--92, Topology and measure 3,1 (1982), Greifswald GDR.*.
The direct appointment transforms
the sequence $\VF\in X^\Phi$ to  a partial set
$\VF\circ \rho \subset X$ which is
called a \normalref{=limit} of the taken sequence.
We shall use the notation
$$\lim_X\ \VF:=\VF\circ \rho\subset X\tsk$$

The points of the limit will be called \normalref{=limit points}.
The opposite appointment transforms the point $x\in X$
 to a set of sequences $\rho\circ x\subset X^\Phi$.
The sequences from this set
will be called \normalref{=converging
to the taken point}.

A \normalref{countable =sequence} in the space $X$ will be a point from the
countable set product of the world space $\VF\in X^\VN$,
i.~e. it is a function
over the natural numbers to the world space
$$\VF_n\in X\uparrow n\in\VN\tsk$$
A natural number $n\in \VN$ will be  called an index
and a value $\VF_n\in X$ will be called  a $n$-member of the sequence.

 The algebraic operations traditionally
 deals with the convergence of \normalref{finite =sequences}
taken from the finite set product $\VF\in X^n$, i.e. a sequence in this
case coincides with a
$n$-tuple of space points $x_i\in X$.

The \normalref{=filter} is a collection of partial sets $\VF\subset X^P$.
The space of all filters will be the second potential of the space
$X^{PP}$. Filters commonly are considered as generalized sequences.

For us any \normalref{=set} may also be considered as a sequence.
The space of such
sequences coincides with potential of the space $X^P$.
The limit point
of such convergence may be interpreted as integral of the taken set,
for example a geometric center point could be  set's limit point
$$x=\int A\tsk$$

For any reform to another world space $R:X\ventrightarrow Y$
we get the reform between the spaces of sequences
$R^\VN:X^\VN\ventrightarrow Y^\VN$ using fractional
extension of mappings. Every mapping $f:X\longrightarrow Y$ provides the
mapping for the sequences $f^\VN:X^\VN\longrightarrow Y^\VN$ coinciding
with set product. For reform we take the fractional extension of set product.
It can be expressed in such way: For every sequence with
members $x_n\in X$ a sequence with members $y_n\in Y$
will be related exactly then, when
these members are related
$$\langle x_n,y_n\rangle \in R\uparrow n\in \VN\tsk$$
The direct appointment appoints to the sequence $x_\cdot\in X^\VN$
 the set coinciding with the product of sets apointed to the members
of the sequence
$$\langle x_n:n\in \VN\rangle \circ R^\VN=\prod\lbrace  x_n\circ R: n\in\VN\rbrace \tsk$$

We describe needed transformations also for another generalized sequences.
For the partial sets $A\in X^P$ we take fractional extension of
mapping's direct image appointment.
 Therefore for a partial set $A\subset X$ we get a related
partial set $B\subset Y$ in another world space  exactly then, when we can
find the partial set $C\subset Q\subset X\times Y$
in the graphic set of the reform $R$ which for projections has direct images
coinciding with the taken partial sets
$$p_1^P(C)=A\kbl p_2^P(B)=B\tsk$$

For the filters $\VF\subset X^{PP}$ we take fractional extension
of mapping's secondary direct image appointment. This secondary image
appointment is used to define the image of the filter for a mapping of the
world space. For a filter $\VF\in X^{PP}$ we have a related
filter in another world space $\VG\in Y^{PP}$ exactly then, when we can find
a filter on the graphic of the reform $\VH\in Q^{PP}$ which
has the images coinciding with the taken filters
$$p_1^{PP}(\VH)=\VF\kbl p_2^{PP}(\VH)=\VG \tsk$$

The continuous function between two convergence spaces
$u: X\longrightarrow Y$ is any function which maintains the limits
$$u^P(\lim_X\ \VF)\subset \lim_Y\  u^\Phi(\VF)\tsk$$

For an algebraic operation, as for function over set product,
the limit of (finite) sequence has exactly
one limit point. Therefore the continuity for such convergence coincides
with equality of limit points. For example, in additive monoid
we demand the  property of additivity
$$f(x+y)=f(x)+f(y)\tsk$$

A map between algebras
is called  \normalref{=morphism of algebras}
if it is continuous for all algebra's  operations.

For example a morphism between  semigroups
$u:X\longrightarrow Y$ is defined
demanding the equality of unique limit points
$$u( \VF\circ \rho_X)=u^\Phi(\VF)\circ \rho_Y \tsk$$
This property is drawn as commuting square
$$\xymatrix{ X^\Phi \ar[r]^-{u^{\Phi}} \ar[d]_{\rho_X} & Y^\Phi \ar[d]^{\rho_Y}
  \\ X \ar[r]_-{u} & Y} $$
We shall say that commuting diagram defines \normalref{=predicate}. In this case
a predicate will be equality between two terms get with
composition of arrows
      $$\rho_X\circ u=u^\Phi\circ \rho_Y\tsk$$

For more general convergence the continuous mappings are defined by
by inclusion of limit sets
$$u^P(\VF\circ\rho_X )\subset u^\Phi(\VF)\circ \rho_Y \tsk$$
Such inclusion can be drawn as \normalref{=directed square}.
It can be named as \normalref{allegory's commuting =square}
$$\xymatrix{
  X^{\Phi} \ar[r] ^(-0.2)*{\gimimas} ^-{u^{\Phi}} \ar[d]_{\rho_X} & Y^{\Phi}
  \ar[d]^{\rho_Y} \\ 2^X \ar[r]_-{u^P} _(1.2)*{\mazanearrow} & 2^Y }$$ 

We can see that this inclusion is equivalent with another inclusion
for opposite reforms
$$ u\inv(y)\circ \rho_X\subset (u^\Phi)\inv(y\circ \rho_Y )\tsk$$
It can be drawn by another directed square
$$\xymatrix{
  X^{\Phi} 
  & Y^{\Phi} \ar[l]_{u^{*\Phi}} _(1.2)*{\mazanearrow} 
  \\ 2^X \ar[u]^-{\rho^*_X} 
  &  2^Y \ar[l]^{u^{*P}}  ^(-0.2)*{\gimimas} \ar[u]_{\rho^*_Y}}$$
We shall say that these two directed squares
define two different \normalref{=reversed predicates}.

\Partnumber,
For diagrams of reforms the predicate is
defined by choosing the \normalref{=birth vertex},
the \normalref{=death vertex},
the two \normalref{=paths from the birth vertex to the death vertex}
and the \normalref{=arrow of inclusion}
between such paths.
The path is \normalref{=cyclic} if it has some repeated vertex.
The path is \normalref{=simple} if it hasn't cyclic parts.
It is a unique possibility for simple path to go from the birth
vertex to the death vertex. Therefore the predicates with simple paths
are determined by choosing the birth and the death vertexes and the
direction of inclusion.

If we denote the vertexes of commuting square
$$\xymatrix @=1.5em{ 1 \ar[r] \ar[d] &2 \ar[d] \\ 4 \ar[r]  &3 }$$
and choose clockwise the positive direction of inclusion, then
the former first directed square will be noted $(-13)$ and the former
second directed square
will be noted $(+31)$. There are 24 predicates with simple paths get from
such square diagram. Every of such predicates can be viewed as a definition
of some continuity property. Therefore clever mathematician
 would wish to
choose the 24 different names for every of such property.
The linguistic problems usually are not interesting for mathematicians,
but there is no other way to understand this new field of multifunctional
algebra. For new theories we need a new efforts of richer language.
My attempt to propose a vocabulary for englishmen can be corrected
by others more prudent wizards  of English language.

\Partnumber,
At first we shall give only four names initiated by two different
definitions of continuous multivalued mappings. Our language will use
the notion of sequences and their limit points.

\Partnumber:
A reform $F: X\ventrightarrow Y$ will be called \normalref{=respecting}
if for every
sequence in the source space
$\VF\in X^\Phi$ we have in the target space the inclusion  of points related
to the limit points in the source space
$$\VF\circ \rho_X\circ F\subset \VF\circ F^\Phi\circ \rho_Y\kbl$$
we shall write
$$F^P(\lim_X \VF)\subset \lim_Y^P (F^\Phi(\VF))\kbl$$
i.e. every point $y\in Y$, related to the limit point $x\in \lim_X \VF$,
must be a limit point for some related sequence $\VF'$.
The $\lim_Y^P$
for some set of filters notes the  direct image appointment.
Such continuity property is noted by predicate $(-13)$.
This property expresses the regularity of limit points in the source
space. The respective reform can be also understood as
distributive operation.
In earlier example the additive morphism is distributive for a
law of addition.

Such property is useful when we want to construct the algebras having two
operations. For example multiplication can be distributive for
the operation of
sum. Taking all possible products of generaters, we get the set of
members which
sums maintains the products. Every product of sums
will be expressed as a sum of products
$$(x+y)\cdot (z+w)=x\cdot z+ y\cdot z+ x\cdot w + y\cdot w\tsk$$
                                                              
\Partnumber,
A reform $F:X\ventrightarrow Y$
will be called \normalref{=creating}  if for every point
in the source space $x\in X$ we have in the target space
the inclusion of sequences $\VF'$ related to convergent
sequences $\VF$ in the source
space
$$x\circ \rho_X^*F^\Phi\subset x\circ F\circ \rho_Y^*\kbl$$
we shall write
$$F^{\Phi P}(\lim_X^*(x))\subset \lim_Y^{*P} (F(x))\kbl$$
i.e. every sequence related to a $x$-convergent sequence must
be $y$-convergent for some point $y\in Y$ related to taken point $x$.
Such property is noted by predicate $(+42)$.
This property expresses the existence of limit points in the target
space.
\Partnumber,
These two properties for mappings coincide  with usual continuity
property: Every limit point in the source space $x\in \lim_X\VF$
must have image $f(x)\in Y$ which is a limit point of the image sequence
$f(\VF)$. Reverse predicate will say that the
$x$-convergent sequence $\VF$ has the image $f(\VF)$
convergent to the image $f(x)\in Y$ of the taken point.

\Partnumber.
The opposite reform $F^*: Y\ventrightarrow X$
provides additional two predicates. We shall give them
new names.
\Partnumber:
A reform $F: X\ventrightarrow Y$ will be called
\normalref{reversely =cautious}
if for every
sequence in the target space
$\VF'\in Y^\Phi$ we have the inclusion of points $x\in X$ having
related limits in the target space
$$\VF'\circ \rho_Y\circ F^*\subset \VF'\circ F^{\Phi*}\circ \rho_X\kbl$$
we shal write
$$(\lim_Y \VF')\circ F^*\subset \lim_X^P(\VF'\circ F^{\Phi*})\kbl$$
i.e. every point in the source space $x\in X$, related to the
limit point $y\in \lim_Y \VF'$ of the taken sequence,
must be a limit point for some sequence $\VF$ in the source
space $X$ for which the taken sequence is related.
Such property is noted by predicate $(+24)$.
Philologically the convergence
in the source space must be cautious as the convergence in the target space
could be respecting.

\Partnumber,
A reform $F: X\ventrightarrow Y$  will be called
\normalref{reversely =wasting}
if for every
point in the target space $y\in Y$ we have inclusion of sequences having
related sequences with taken point in the target space
$$y\circ \rho_Y^*\circ F^{\Phi*}\circ \rho_X^* \kbl$$
we shall write
$$(\lim_Y^* (y))\circ F^{\Phi*}\subset \lim_X^{*P}(y\circ F^*)\kbl$$
i.e. every sequence $\VF$ in the source space having
related sequence $\VF'$ in the
target space with taken limit point $y\in\lim_Y \VF'$, must be convergent
to some point $x\in X$ for which taken point $y\in Y$ is related.
Such property is noted by predicate $(-31)$.
Philologically the convergence in the source space must be wasting, as
the convergence
in the target space could be creating.

\Partnumber.
The opposition maintains the inclusion of the graphics, therefore
we can define additionally 4 predicates with known properties.
The predicate $(+31)$ notes the \normalref{reversely =respecting}
reform $F$,
the predicate $(-24)$ notes the \normalref{reversely =creating} reform $F$,
the predicate $(-42)$ notes the \normalref{=cautious} reform $F$,
and the predicate $(+13)$ notes the \normalref{=wasting} reform $F$.
These predicates will define different properties for pairs of conjugate
functors. For reforms such functors are given by direct image appointment
and opposite image appointment. Therefore we wished to have different names
for all predicates. At this moment I choose  the
predicate with direct reform as principal.

\Partnumber,
We have defined the names for 8 predicates. We shall say that
these predicates are from the \normalref{first =octet}. When the reform
$F:\ventrightarrow Y$ is a mapping, we can drive along this mapping
the death or the birth
vertexes for the first two predicates $(-13)$, $(-24)$. We get
the new predicates $(-14)$ and $(-23)$, which defines the same continuity
properties for mappings.
These predicates begins the \normalref{second =octet}.
\Partnumber:
A reform $F: X\ventrightarrow Y$ will be called \normalref{=lavishing} if
for every sequence $\VF$ in a source space $X$ every limit point has related
point in the target space $y\in Y$ which is a limit point for some related
sequence $\VF'$
$$\VF\circ \rho_X\subset \VF\circ F^\Phi\circ \rho_Y\circ F^*\kbl$$
we can write
$$\lim_X\VF\subset (\lim_Y^P (\VF\circ F^\Phi))\circ F^*\tsk$$
It is more convenient to express such inclusion as large opposite image
of convergence in the target space
$$\rho_X\subset (F^\Phi\times F)^{*P}(\rho_Y ) \tsk$$
For the reforms there we apply the fractional extension of set product.
Such property is noted by predicate $(-14)$.

\Partnumber,
A reform $F: X\ventrightarrow Y$ will be called \normalref{=pressing} if
for every sequence $\VF'$ in a target space $Y$ the limit point $x\in X$ of
sequences $\VF$ in a source space,
 which has related taken sequence $\VF'$,
has related some limit point $y\in \lim_Y \VF'$
$$\VF'\circ F^{\Phi*}\circ \rho_X\circ F\subset \VF'\circ \rho_Y\kbl$$
we shall write
$$(\lim_X^P (\VF'\circ F^{\Phi*}))\circ F\subset \lim_Y\VF'\tsk$$
It can also be expressed that a  convergence
in a source space has  small direct image
$$(F^\Phi\times F )^P(\rho_X )\subset \rho_Y\tsk$$
Such property is noted by predicate $(-23)$.

\Partnumber.
The opposite reform $F^*: Y\ventrightarrow X$ provides a new names.
\Partnumber:
A reform $F: X\ventrightarrow Y$ will be called \normalref{=hiding}
if for every
sequence $\VF'$ in the target space $Y$ every limit point $y\in Y$ is related
to some point in the source space $x\in X$, which is a limit point for
some sequence $\VF$ in a source space, which has the taken sequence $\VF'$
as related
$$\VF'\circ \rho_Y\subset \VF'\circ F^{\Phi *}\circ \rho_X\circ F\kbl$$
we shall write
$$\lim_Y \VF'\subset (\lim_X^P (\VF'\circ F^{\Phi*}))\circ F\tsk$$
It can be expressed that the convergence in a source space
has large direct image
$$\rho_Y\subset (F^\Phi\times F)^P(\rho_X)\tsk$$
Such property is noted by predicate $(+23)$.
Philologically the convergence in a source space must be hiding, as the
convergence in a target space could be lavishing.

\Partnumber,
A reform $F: X\ventrightarrow Y$ will be called \normalref{=reflecting}
if for
every sequence $\VF$ in the source space $X$ every point $x\in X$,
which has related limit point of some related sequence $\VF'$ in a target
space, is a limit point of taken sequence $\VF$
$$\VF\circ F^\Phi\circ \rho_Y\circ F^*\subset \VF\circ \rho_X\kbl$$
we shall write 
$$(\lim_Y^P (\VF\circ F^\Phi))\circ F^*\subset \lim_X\VF\tsk$$
It can be expressed that the convergence in a target space has a
small opposite image
$$(F^\phi\times F)^{*P}(\rho_Y)\subset \rho_X\tsk$$
Tis property is noted by predicate $(+14)$.

\Partnumber.
The reversing will give the same names for additional 4 predicates.
At this moment the principal predicate deals with the limit points of
convergence, and
reversed predicate deals with convergent sequences of convergence.
The predicate $(+41)$ notes the \normalref{reversely =lavishing} reform $F$,
the predicate $(+32)$ notes the \normalref{reversely =pressing} reform $F$,
the predicate $(-32)$ notes the \normalref{reversely =hiding} reform $F$,
the predicate $(-41)$ notes the \normalref{opposite =reflecting} reform $F$.

\Partnumber,
We have just defined the 8 predicates from the second octet.
The rest predicates compound the \normalref{third =octet}.
For functional convergence we can get the earlier continuity properties.
This helps us to choose the preference order in description of
last predicates.
\Partnumber:
A reform $F: X\ventrightarrow Y$ will be called \normalref{=thin}
if for every sequence $\VF$ in a source space $X$ every sequence converging
to a point, related to some limit point $x\in X$ of taken sequence $\VF$,
is related with taken sequence $\VF$
$$\VF\circ \rho_X\circ F\circ \rho_Y^*\subset \VF\circ F^\Phi\kbl$$
we shall write
$$\lim_Y^{*P} (F^P(\lim_X \VF))\subset F^{\Phi P}(\VF)\tsk$$
It is better to say that the inverse image of the taken reform is  smaller
than the reform for the sequences
$$(\rho_X\times \rho_Y)^{*P}(F)\subset F^\Phi\tsk$$
This property is noted by predicate $(-12)$.
\Partnumber,
A reform $F: X\ventrightarrow Y$ will be called \normalref{=thick}
if for every point $x\in X$ in a source space the limit points in target
space $y\in\lim_Y \VF'$, of some sequence related to some sequence $\VF$
in source space converging to the taken point $x$,
is related to the taken point $x$
$$x\circ \rho_X^*\circ F^\Phi\circ \rho_Y\subset x\circ F \kbl$$
we shall write
$$\lim_Y^P( F^{\Phi P}(\lim_X^*(x)) \subset F(x)\tsk$$
It is better to say that the direct image of the  reform for sequence
is smaller than the taken reform itself
$$(\rho_X\times \rho_Y)^P(F^\Phi)\subset F\tsk$$
This property is noted by predicate $(+43)$.
This property speaks about regularity of existing limit points. Such
reforms were often encountered in the theory of differential operators,
they are called \normalref{ operator with =closed graphic}.
\Partnumber.
The reversing provides additional two names.
The predicate $(+21)$ notes the \normalref{reversely =thin} reform $F$,
and the predicate $(-34)$ notes the \normalref{reversely =thick}  reform $F$.
These predicates deal with graphic of opposite reform $R^*$.
\Partnumber,
Rest 4 last predicates.
\Partnumber:
A reform $F: X\ventrightarrow Y$ will be called \normalref{=binding}
if for every sequence $\VF$ in the source space $X$ every related
sequence $\VF$ in the target space $Y$ has limit point $y\in Y$ related
to some limit point $x\in X$ of taken sequence $\VF$
$$\VF\circ F^\Phi\subset \VF\circ\rho_X\circ F\circ \rho_Y^*\kbl$$
we shall write
$$F^{\Phi }(\VF)\subset \lim_Y^{*P} (F^P(\lim_X \VF ))\tsk$$
It is better to say that the inverse image of the taken reform is larger than
the reform for sequences
$$F^\Phi\subset (\rho_X\times \rho_Y)^{*P}(F) \tsk$$
This property is noted by the predicate $(+12)$.
\Partnumber,
A reform $F: X\ventrightarrow Y$ will be called \normalref{=parting}
if for every point $x\in X$ in a source space $X$ the related point $y\in Y$
in the target space $Y$ must be a limit point of some sequence $\VF$ which
is related to some sequence $\VF$ in a source space $X$ converging
to taken point $x$
$$x\circ F\subset x\circ \rho_X^*\circ F^\Phi\circ \rho_Y\kbl$$
we shall write
$$F(x)\subset \lim_Y^P F^{\Phi}(\lim_X^*(x))\tsk$$
It is better to say that direct image of the sequence reform
is larger than the taken reform itself
$$F\subset (\rho_X\times \rho_Y)^P(F^\Phi)\tsk$$
This property  is noted by the predicate $(-43)$.
\Partnumber.
The reversing provides additional two names.
The predicate $(-21)$ notes the \normalref{reversely =binding} reform $F$,
and the predicate $(+34)$ notes the \normalref{reversely =parting}
reform $F$.
These predicates deals with graphic of opposite reform $R^*$.

\Partnumber,
The reform composition maintains the inclusion
$$F\subset G\duoda R\circ F\subset R\circ G, F\circ H\subset G\circ H\tsk$$
This property gives a
possibility to apply the composition of directed squarees. This is
straightforward applied for the predicates of first octet.

\Partnumber:
\Proposition.
The composition of reform maintains the property of
respecting reforms $(-13)$,
creating reforms $(+42)$, cautious reforms $(+24)$,
reversely wasting reforms $(+31)$ and their reversed counterparts.
\Proof:        
It is enough to apply composition of \normalref{=neighbouring squares}.
\irodymopabaiga
\Partnumber,
The predicates of second octet needs the composition of
\normalref{=embracing squares}.
\Proposition.
The composition of reform maintains
the property of lavishing reforms $(-14)$, pressing reforms $(-23)$,
hiding reforms $(+23)$, reflecting reforms $(+14)$ and
their reversed counterparts.
\Proof:
We construct the embracing squarees.
\irodymopabaiga
\Partnumber,
For the predicates of third octet we must pose restrictions for the
convergence in the mediating space. At first we explain the terminology.

For every reform $R:X\ventrightarrow Y$ we have defined two
appointments: the direct image $R^P: X^P\longrightarrow Y^P$
and inverse image $R^{*P}: Y^P\longrightarrow X^P$. They are monotonic
mappings of Boole algebras, therefore we have instances of functors.

We shall call the pair of opposite functors $F$ and $G$
between two categories
an \normalref{=adjunction} with the same notation as for reform
$R:\VA\ventrightarrow \VB$. Therefore for this adjunction  we have the
source category $\VA$ and the target category $\VB$. The direct functor
$F:\VA\Longrightarrow \VB$ is called \normalref{=adjoint functor},
and opposite
functor $G:\VB\Longrightarrow \VA$ is called \normalref{=coadjoint functor}.
The \normalref{=unit} of adjunction is a transformation of Identity functor
$i:\Identity\longrightarrow FG$ and \normalref{=counit} of adjunction is a
transformation of functor composition $e: GF\longrightarrow \Identity$.
The naturality and triangular identities is needed for the usual
\normalref{true =adjunction},
but we can deal with ``generalized'' adjunctions without such properties.

The entire reform $R$ provides a unit of its adjunction
$$A\subset R^{*P}(R^P(A))\uparrow A\subset X\kbl$$
the accurate reform $R$ provides a counit of its adjunction
$$R^P(R^{*P}(A'))\subset A'\uparrow A'\subset Y\kbl$$
the covering reform $R$ provides a unit of opposite adjunction
$$A'\subset R^P(R^{*P}(A'))\uparrow A'\subset Y\kbl$$
and the exact reform $R$ provides counit of opposite adjunction
$$R^{*P}(R^P(A))\subset A\uparrow A\subset X\tsk$$

The mapping $R:=f$ provides a unit and counit together, therefore such
adjunction for Boole algebras is true. The inverse of
mapping $R:=f^*$ provides a unit and
counit for opposite adjunction, therefore in this case the opposite
adjunction is true. For bijective map $R$ both adjunctions are true.

\Proposition.
The composition maintains the property of thin reforms $(-12)$
if the convergence in the middle space is simple.
The composition maintains the property of thick reforms $(+43)$
if the convergence in the middle space is entire.
The composition maintains the property of binding reforms $(+12)$
if the convergence in the middle space is covering.
The composition maintains the property of parting reforms $(-43)$
if the convergence in the middle space is exact.
The reversed counterparts needs the same conditions for the convergence
in the middle space.
\Proof:
These conditions provide needed unit and co-unit arrows of direct
and opposite adjunctions.
\irodymopabaiga
For the countable sequence $x_n\in X\uparrow n\in X$
an accurate  convergence is known as Hausdorff.
Usual sequence covergence is covering, as it appoints the point
$x\in X$ for constant
sequence $x_n=x$. The usual  convergence is not exact, because different
sequences can converge to the same limit point.
Entire convergence may be only for some algebraic operation over the
finite part of sequences.

\Partnumber.
The identity mapping has continuity property of every kind, therefore
the respective reforms compound the suballegory. Sometimes it is beneficial
to use the decomposition
with tabulating mappings $R=f\inv\circ g$ in this suballegory, because
the continuity properties for mappings is much simple.

The death vertex can be pushed forward along adjoint functor with counit, and
backword along adjoint functor with unit. The direction is defined
with the arrow of predicate. In the case of set reforms we get
that death vertex can be pushed forward along the accurate reform, and
can be pushed backward along the entire reform.

The birth vertex can be pushed forward along coadjoint functor with
unit, and can be pushed backward along coadjoint functor with counit.
In the case of set reforms we get that the birth vertex can be pushed
forward along covering reform, and can be pushed backward along
exact reform.

After pushing we get the continuity for the new predicate.
In the case of mappings we get the equivalent continuity properties.
The death vertex can be pushed forward and backward along mappings, and
birth vertex can be pushed backward and forward along inverse mappings.

Therefore we can formulate the proposition about continuity properties
of mappings.

\Proposition.
For mappings 4 continuity properties are equivalent: respecting $(-13)$,
pressing $(-23)$, lavishing $(-14)$ and reversely creating $(-24)$.
\Proof:
Such predicates allow to push the death and birth vertexes along the
mappings or inverse mappings.
\irodymopabaiga
This proposition also is valid for 4 rversed  predicates.


%% file: madher.tex
\Part Reform's continuity in adherence spaces.
\Partnumber:
The \normalref{=adherence} in the space $X$ is defined as generalized
convergence
$\Adh:2^X\ventrightarrow X$. A sequence now will be an arbitrary partial set
$A\subset X$. The direct appointment appoints the set of \normalref{=adherent points}
$$[A]:=A\circ \Adh$$
and will be called an \normalref{=adherence operator}.
The opposite appointment  appoints the filter of sets having
taken adherent point
$$ \VF_x:=\Adh\circ x\tsk$$
We shall call such filter an \normalref{=adherence filter over the point}
 $x\in X$.

We remember that any of these appointments restore the adherence.
We shall say the point $x\in X$ \normalref{is =adherent}
to the set $A\subset X$
and the set \normalref{=adherent point}.
All adherent points for the set $A\subset X$
compound the \normalref{=adherence set} $[A]\subset A$.

The \normalref{=conjugate operator} is called \normalref{=inside operator},
it is defined with
the Boole completement of the sets
$$\langle A\rangle :=A\circ \Ins= (A^c\circ \Adh)^c \tsk$$
The sets with the adherent point now will be changed by the neighbourhoods
of taken point. The set $V\subset X$ is a neighbourhood of the point $x\in X$
and this point $x\in X$ is inside  point of this set
exactly then, when the completment
of this set $A^c=X\setminus A$ hasn't taken point as adherent point.
We shall write $$\VV_x:=\Ins\circ x\tsk$$
We shall say that the inside points $x\in A$ compound the inside
set $\langle A\rangle \subset X$. The inside operator or
neighbourhood filters
define \normalref{=inside} as  adherence conjugate to taken one.

It is obvious that secondly conjugate operator coincides with the taken one.
The names adherence or inside are only relative consideration. At this
moment every adherence can be inside of the dual adherence. Next we can
ask some asymmetric properties for the adherence.

\Partnumber:
An \normalref{=increasing adherence} appoints the bigger
adherence set $A\subset [A]$.
An \normalref{=isotonic adherence} for the bigger set apoints a bigger set
$$A\subset B\duoda [A]\subset [B]\tsk$$
The isotonic increasing adherence is called a \normalref{=closure}.
The set is called closed
if it coinsides with it's adherence set
$$A=[A]\tsk$$
For such adherence the inside will be \normalref{=diminishing}
$A\supset \langle A\rangle $
and isotonic
$$A\subset B\duoda \langle A\rangle \subset \langle B\rangle \tsk$$
The isotonic diminishing inside will be called \normalref{=interior}.
The set
is open if it coincides with it's interior
$$A=\langle A\rangle \tsk$$

\Partnumber.
We shall look the continuity properties for the reforms of
adherence or inside spaces.
We shall content himself with the case of mappings
$f: X\longrightarrow Y$.
Then the mapping of partial sets is provided by direct image
$f^P:X^P\longrightarrow Y^P$. We shall use the predicates from the first
and second octets only. Such predicates are good for the sequential
adherences, but we haven't opportunity to speak there about this subject
now. The adherence operator and the adherence filters helps us to
deal continuity properties more easy.

At first we shall check the equivalent
predicates which coincide with property of traditional continuous
mapping.
\xymatrixrowsep={1.6pc}
\xymatrixcolsep={2.6pc}
\Partnumber:
The respecting  $(-13)$ mapping of adherence spaces is defined by directed
squaree:
$$\xymatrix {\centre{-13}
  \gimimas \ar[r]^-{f^P} \vent[d]_{\Adh_X} 
  & \cdot \vent[d] ^{\Adh_Y}
  \\ \cdot  \ar[r]_-{f} _(1.2)*{\mazanearrow} & \cdot }$$
The corresponding inclusion of reforms
$$A\circ \Adh_X\circ f\subset A\circ f^P\circ \Adh_Y $$
is easily defined  with adherence operators
i.e. we get usual condition for continuous mappings in topological
spaces
$$f^P([A])\subset [f^P(A)]\uparrow A\subset X\tsk $$

The reversely  respecting $(+31)$ mapping is easily  defined
with adherence filters:
$$\xymatrix{\centre{+31}
  \cdot 
  & \cdot \vent[l]_{f^{P*}} _(1.2)*{\mazanearrow} 
  \\ \cdot \vent[u]^{\Adh^*_X}  
  & \gimimas  \vent[l]^-{f^*} \vent[u]_{\Adh^*_Y}}$$
The inclusion of reforms
$$y\circ f^*\circ \Adh_X^*\subset y\circ \Adh_Y^* \circ f^{P*}$$
can be written as properties of adherence filters
$$\VF_{f^{-1}(y)}\mzl (f^P)\inv (\VF_y)\tsk$$
For every point in a target space $y\in Y$ every point of inverse image
$x\in f\inv(y)$ is adherent only for such set $A\subset X$
which has direct image $f^P(A)\subset Y$ for which taken point is adherent
$y\in [f^P(A)]$.
\Partnumber,
The reversely creating  $(-24)$ mapping of adherence spaces
$$\xymatrix{\centre{-24}
  \cdot  \vent[d]_{\Adh_X} 
  & \gimimas \vent[l]_{f^{P*}} \vent[d]^{\Adh_Y} 
  \\ \cdot  & \cdot \vent[l]^-{f^*} ^(1.2)*{\mazasearrow} }$$
is easily defined with adherence operators.
The inclusion of reforms
$$B\circ f^{P*}\circ \Adh_X\subset B\circ \Adh_Y\circ f^*$$
we want to change with inclusion of adherence sets. For every partial
set in the source space $A\subset X$ having the direct image coinciding
with taken partial set in the target space $B=f^P(A)$, the adherence
set is smaller than inverse image of taken partial set
$$A\subset f\inv(B)\tsk$$
For isotonic adherences this property can be get more easily.

\Proposition.
For reversely creating mapping $f:X\longrightarrow Y$
every covered partial set $B\subset f^P(X)$ has inclusion for
adherence set of inverse image 
$$[f\inv (B)\subset f\inv ([B])\tsk$$
If the adherence in the source space is isotone, then this property
is also sufficient for mapping to be reversely creating.

If the adherence in the target space is isotone, then this property
is equivalent to the same requirement for every partial set in the target
space $B\subset Y$.
\Proof:
For covered set $B\subset f^P(X)$ the inverse image $A=f\inv(B)$ will have
the image coinciding with the taken set $f^P(A)=B$, therefore its adherence
set is included in the inverse image of adherence set for the taken set
$$[f\inv (B)]\subset f\inv([B])\tsk$$

For every partial set in the target space $B\subset Y$ we can take smaller
set $B'=f^P(f\inv (B))$, and for isotone adherence in the target space
$$[f\inv (B)]=[f\inv (B')]\subset f\inv ([B'])\subset f\inv ([B])\tsk$$

For every partial set $A\subset X$ with image $f^P(A)=B$
we have $A\subset f\inv (A)$, therefore for isotone adherence in the source
space we shall have
$$[A]\subset [f\inv(B)]\subset f\inv ([B])\tsk$$
\irodymopabaiga

For isotonic adherences this property is equivalent with the inclusion
of adherence for the inverse image
$$[f\inv (B)]\subset f\inv ([B])\tsk$$

The creating $(+42)$ mapping
$$\xymatrix{\centre{+42}
   \cdot 
   & \cdot \ar[l]_{f^P} _(-0.2)*{\mazasearrow}  
   \\ \gimimas \vent[u]^{\Adh_X^*} 
   & \cdot \vent[u]_{\Adh_Y^*} \ar[l]^-{f} }$$
is defined by inclusion  of reforms
$$x\circ \Adh_X^*\circ f^P\subset x\circ f\circ \Adh_Y^*\tsk$$
This can be written with inclusion of image for adherence filter
$$f^{PP}(\VF_x)\mzl \VF_{f(x)}\tsk$$
Every partial set $A\subset X$ adherent by the taken point $x\in X$
has image $f^P(A)\subset Y$ adherent by the image of the taken point
$y=f(x)$.

\Partnumber,
The pressing $(-23)$ mapping
$$\xymatrix{\centre{-23}
  \cdot \vent[d]_{\Adh_X} 
  & \gimimas \vent[l]_{f^{P*}} \vent[d]^{\Adh_Y}
  \\ \cdot \ar[r]_-{f} _(1.2)*{\mazanearrow} & \cdot }$$
is defined by inclusion of
reforms
$$B\circ f^{P*}\circ\Adh_X\circ f\subset B\circ \Adh_Y\tsk$$
It can be expressed as small image of adherence $\Adh_X\subset X^P\times X$
$$(f^P\times f)^P(\Adh_X)\subset \Adh_Y\tsk$$
For the adherence operator we shall get small images
$$f^P([A])\subset [B]$$
for the adherent set of
sets $A\subset X$ which has direct image
coinciding with taken set $f^P(A)=B$. This  can be written as
inclusion of the filters
$$f^{PP}([(f^P)\inv (B)])\mzl \lbrace [B]\rbrace \tsk$$
For isotone adherences this property is expressed more easily.

\Proposition.
For the isotone adherence in the source space the mapping is pressing
if every covered partial set $B\subset Y$ has
small image of adherence set for the inverse image of the taken set
$$f^P([f\inv(B)])\subset [B]\tsk$$
For the isotone adherence in the target space this requirement is
equivalent to such inclusion for every partial set $B\subset Y$.
\Proof:
For isotone adherence in the source space a smaller partial set will
have a smaller adherence set
$$A\subset f\inv(B)\duoda [A]\subset [f\inv (B)]\kbl$$
therefore we shall have the inclusion for the image of adherence set
$$f^P([A])\subset f^P([f\inv(B)])\tsk$$
We have shown that for such adherence every partial set $A\subset X$
will have a small image of adherence set
$$f^P([A])\subset [B]\tsk$$
If the adherence in the target space is isotone, then for every partial set
$B\subset Y$ we have a smaller covered set
$f^P(f\inv (B))\subset B$,
and we can check the needed inclusion
$$f^P([f\inv (B)])\subset [f^P(f\inv (B))]\subset [B ]\tsk$$
\irodymopabaiga
We have shown
that
for isotone adherences the property of pressing mapping
is expressed 
with inverse image of taken set in the target space $$f^P([f\inv (B)])
\subset [B]\uparrow B\subset Y\tsk$$

The reversely pressing $(+32)$ mapping
$$\xymatrix{\centre{+32}
  \cdot \ar[r]^-{f^P} ^(1.2)*{\mazasearrow}
  & \cdot 
  \\ \cdot \vent[u]^{\Adh_X^*} 
  & \gimimas \vent[l]^-{f^*} \vent[u]_{\Adh_Y^*}}$$
is defined by inclusion of opposite reforms
$$y\circ f^*\circ\Adh_X^*\circ f^P\subset y\circ \Adh_Y^*\kbl$$
which also express small image of adherence
$$(f\times f^P)^P(\Adh_X^*)\mzl \Adh_Y^*\tsk$$
For adherence filters we get that for any point in a target space $y\in Y$
the points of inverse image $x\in f\inv (y)$ have adherence filter $\VF_x$
with small image
$$f^{PP}(\VF_{f\inv(y)})\mzl \VF_y\tsk$$

\Partnumber,
The lavishing $(-14)$ mapping
$$\xymatrix{\centre{-14}
  \gimimas \vent[d]_{\Adh_X} 
  & \cdot \ar[l]_{f^P} \vent[d]^{\Adh_Y}
  \\ \cdot & \cdot \vent[l]^-{f^*} ^(1.2)*{\mazasearrow}}$$
demands an inclusion of reforms
$$A\circ \Adh_X\subset A\circ f^P\circ \Adh_Y\circ f^*\kbl$$
i.e. we ask a large inverce image of the adherence
 $$\Adh_X\subset (f^P\times f)\inv (\Adh_Y)\tsk$$
For adherence operator it will be  an inclusion of adherence set
$$[A]\subset f\inv ([f^P(A)])\uparrow A\subset X\tsk$$

The reversely lavishing $(+41)$ mapping
$$\xymatrix{\centre{+41}
  \cdot 
  & \cdot \vent[l]_{f^{P*}} _(1.2)*{\mazanearrow}
  \\ \gimimas \vent[u]^{\Adh_X^*} \ar[r]_-{f} 
  & \cdot \vent[u]_{\Adh_Y^*} }$$
asks another inclusion
$$\Adh_X^*\subset f\circ \Adh_Y^*\circ f^{P*}\kbl$$
which again is equivalent of asking the large inverse image of adherence
$$\Adh_X^*\subset (f\times f^P)\inv(\Adh_Y^*)\tsk$$
It can be expressed with large inverse image of adherence filter
$$\VF_x\mzl (f^P)\inv (\VF_{f(x)})\tsk$$

\Partnumber.
Remain the predicates without equivalent counterparts.
\Partnumber:
Wasting $(+13)$ mapping
$$\xymatrix{\centre{+13}
  \gimimas \vent[d]_{\Adh_X} \ar[r]^-{f^P} 
  & \cdot \vent[d]^{\Adh_Y}  
  \\ \cdot \ar[r]_-{f} _(1.2)*{\mazaswarrow} & \cdot }$$
asks an inclusion of reforms
$$A\circ f^P\circ \Adh_Y\subset A\circ \Adh_X\circ f \tsk$$
It can be expressed with adherence operators
$$[f^P(A)]\subset f^P([A])\uparrow A\subset X\tsk$$

The reversely wasting $(-31)$ mapping
$$\xymatrix{\centre{-31}
  \cdot 
  & \cdot \vent[l]_{f^{P*}} _(1.2)*{\mazaswarrow}
  \\ \cdot \vent[u]^{\Adh_X^*} 
  & \gimimas \vent[l]^-{f^*} \vent[u]_{\Adh_Y^*}}$$
asks an inclusion of reforms
$$y\circ\Adh_Y^*\circ f^{P*}\subset y\circ f^*\circ\Adh_X^*\tsk$$
This can be expressed with small inverse image of adherence filter
$$(f^P)\inv(\VF_y)\subset \VF_{f\inv(y)}\tsk$$
\Partnumber,
The reversely cautious $(+24)$ mapping
$$\xymatrix{\centre{+24}
  \cdot \vent[d]_{\Adh_X} 
  & \gimimas \vent[l]_{f^{P*}} \vent[d]^{\Adh_Y}
  \\ \cdot & \cdot \vent[l]^-{f^*} ^(1.2)*{\mazasearrow}}$$
is asking the inclusion  of reforms
$$B\circ \Adh_Y\circ f^*\subset B\circ f^{P*}\circ \Adh_X\tsk$$
This property can be expressed with adherence operator:
Every point $x\in f\inv([B])$ in the inverse image of the adherence set
for taken set in the target
space $B\subset Y$ is adherent $x\in [A]$ for some partial set in
the source space
$A\subset X$  with image coinsiding with taken set $f^P(A)=B$.

\Proposition.
For the isotone adherence in a source space the mapping
$f:X\longrightarrow Y$ will be reversely  cautious exactly then,
 when for every
covered partial set $B\subset f^P(X)$ the inverse image of the adherence
set would contained in the adherence set of the inverse image of
taken set
$$f\inv([B])\subset [f\inv (B)]\kbl$$
and for other partial sets in a target space $B\not \subset f^P(X)$
the adherence set would contained in completment of mapping's whole image
$$[B]\subset Y\setminus f^P(X)\tsk$$
\Proof:
This condition is obviously sufficient.

And it is necessary:
If a point $x\in X$ is in inverse image of the adherence set
$x\in f\inv([B])$ for a
covered set $B\subset f^P(X)$, then any set $A\subset X$ with the image
$f^P(A)=B$ and adherent by the taken point $x\in [A]$ will have
$A\subset f\inv (B)$, therefore for isotone adherence
$$x\in [A]\subset [f\inv (B)]\tsk$$

For other partial sets $B\not\subset f^P(X)$ the inverse image of the
adherence set will be empty
$f\inv ([B])=\tustis$, as we can't find any partial set $A\subset X$
having image equal to the taken set, i.e.
$$f^P(A)\not=B\uparrow A\subset X\tsk$$
\irodymopabaiga

The cautious $(-42)$ mapping
$$\xymatrix{\centre{-42}
  \cdot & \cdot \ar[l]_{f^P} _(-0,2)*{\mazanwarrow}
  \\ \gimimas \ar[r]_-{f} \vent[u]^{\Adh_X^*} 
  & \cdot \vent[u]_{\Adh_Y^*} }$$
asks an inclusion of reforms
$$x\circ f\circ \Adh_Y^*\subset x\circ\Adh_X^*\circ f^P\tsk$$
This is expressed with big image of the adherence filter
$$f^{PP}(\VF_x)\dgl \VF_{f(x)}\tsk$$
It is very special property for the mapping: the point's image $f(x)\in Y$
can be adherent only such partial sets $B\subset Y$ which coincide
with direct image of any partial set in the source space
$B=f^P(A)$.

\Partnumber,
The hiding $(+23)$ mapping
$$\xymatrix{\centre{+23}
  \cdot \vent[d]_{\Adh_X}
  & \gimimas \vent[l]_{f^{P*}} \vent[d]^{\Adh_Y}
  \\ \cdot \ar[r]_-{f} _(1.2)*{\mazaswarrow} 
  & \cdot }$$
demands an inclusion of relations
$$B\circ \Adh_Y\subset B\circ f^{P*}\circ \Adh_X\circ f\kbl$$
which equivalent to the large image of adherence
$$\Adh_Y\subset (f^P\times f)^P(\Adh_X)\tsk$$
This can be expressed  with adherence operator. For the adherent point
$y\in [B]$ we can find partial set in a source space $A\subset X$
which image coincides with the taken set
$f^P(A)=B$ and we have adherent point $x\in [A]$ covering the earlier point
in the target space $f(x)=y$
$$[B]\subset \Cup\lbrace f^P([A]): f^P(A)=B\rbrace \tsk$$

For the isotone adherence in the source space the hiding mapping
can be defined by stronger requirement.
\Proposition.
If the whole target set has any adherent point $[Y]\not=\tustis$,
then hiding mapping $f:X\longrightarrow Y$ will cover the whole target
space.

For the isotone adherence in the source space $\Adh_X$
the hiding mapping $f:X\longrightarrow Y$ will
provide every covered partial set $B\subset f^P(X)$ with adherence set
contained in the image of adherence set for inverse image of taken set
$$[B]=\tustis\uparrow B\not \subset f^P(X)\tsk$$
\Proof:
For not covered set $B\subset Y$ we cant find any partial set in the source
space $A\subset X$ with coinciding image $f^P(A)=B$, therefore
for hiding mapping such set cannot have any adherent point $[B]=\tustis$.

For the covered set $B\subset Y$ every adherent point $y\in [B]$
must be covered by adherent point $x\in [A]$ for some partial set in the
source space $A\subset X$ with coinciding image $f^P(A)=B$.
If the adherence in the source space is isotone, then
$$y\in f^P([A])\subset f^P([f\inv (B)])\uparrow B\subset Y\tsk$$
\irodymopabaiga

The reversely hiding $(-32)$ mapping
$$\xymatrix{\centre{-32}
  \cdot \ar[r]^-{f^P} ^(1.2)*{\mazanwarrow}
  & \cdot  
  \\ \cdot \vent[u]^{\Adh_X^*} 
  & \gimimas \vent[u]_{\Adh_Y^*} \vent[l]^-{f^*}}$$
demands an inclusion of reforms
$$y\circ\Adh_Y^*\subset y\circ f^*\circ \Adh_X^*\circ f^P\kbl$$
which is equivalent for large image of adherence's opposite reform
$$\Adh_Y^* \mzl (f\times f^P)^P(\Adh_X^*)\tsk$$
This will be expressed as large image  of adherence filter
appointment
$$\VF_y\mzl \Cup \lbrace f^{PP}(\VF_x):x\in f\inv(x)\rbrace \tsk$$

We have seen that such mapping must be very special.
For the target space with
adherent point $[Y]\not=\tustis$ such mapping must be covering.

\Partnumber,
The reflecting $(+14)$ mapping
$$\xymatrix{\centre{+14}
  \gimimas \ar[r]^-{f^P} \vent[d]_{\Adh_X}
  & \cdot \vent[d]^{\Adh_Y}
  \\ \cdot 
  & \cdot \vent[l]^-{f^*} ^(1.2)*{\mazanwarrow}}$$
is defined by inclusion of reforms
$$A\circ f^P\circ \Adh_Y\circ f^*\subset A\circ \Adh_X \kbl$$
which is equivalent to the small inverse image of the adherence's graphic
$$(f^P\times f)\inv (\Adh_Y)\subset \Adh_X\tsk$$
This can be expressed with the adherence operator
as large adherence set in the source space
$$f\inv ([f^P(A)])\subset [A]\uparrow A\subset X\tsk$$
The adherent point of some set $x\in [A]$ is a property which is reflected
by such mapping. If this property is checked in the target space
for image point and the set $f(x)\in [f^P(A)]$, then we can conclude
that this property is also content for initial point and partial set in
the source space.

The reversely reflecting $(-41)$ mapping
$$\xymatrix{\centre{-41}
  \cdot 
  & \cdot \vent[l]_{f^{P*}} _(1.2)*{\mazaswarrow}
  \\ \gimimas \ar[r]_-{f} \vent[u]^{\Adh_X^*} 
  & \cdot \vent[u]_{\Adh_Y^*} } $$
is defined by inclusion of reforms
$$x\circ f\circ \Adh_Y^*\circ f^{P+}\subset x\circ \Adh_X^*\kbl$$
which is equivalent to the \normalref{=small inverse image}
of adherence's opposite
reform
$$(f\times f^P)\inv (\Adh_Y^*)\mzl \Adh_X^*\tsk$$
This property can be expressed for the appointment of adherence filters
$$(f^P)\inv (\VF_{f(x)})\mzl \VF_x\tsk$$

\Partnumber.
Finely we formulate the predicates from the third octet. These
properties is not useful for the traditional topological adherences.
But they may be interesting for some integral representation
of set's middle point. The adherence point $x\in [X]$
can be understood as middle point of this set. Such point usually exists
and is defined uniquely, but we can questionize such existence or
may exist many of needed middle points for some extraordinary sets.

\Partnumber:
The binding $(+12)$ mapping
$$\xymatrix{\centre{+12}
  \gimimas \ar[r]^-{f^P} ^(1.2)*{\mazasearrow} \vent[d]_{\Adh_X}
  & \cdot 
  \\ \cdot \ar[r]_-{f} 
  & \cdot \vent[u]_{\Adh_Y^*} } $$
is defined by inclusion of reforms
$$A\circ f^P\subset A\circ \Adh_X\circ f\circ \Adh_Y^* \kbl$$
which is equivalent to the large inverse image of the taken mapping
$$f^P\subset (\Adh_X\times Adj_Y)\inv (f)\tsk$$
It may be defined with the union  of adherence filters
$$f^P(A)\in  \VF_{f^P([A])}\kbl$$
the image of the partial set $A\subset X$ has the image $f(y)\in Y$ of the
middle point $x\in [A]$ as a middle point.

For functional convergence in the source space
this property would be the same
as for the creating mapping $(+42)$,
and for functional convergence in the target space
this property would be the same as  for the wasting mapping $(+13)$,
therefore it is a same as usual continuity property for mappings
defined with reversely respecting $(+31)$, reversely pressing $(+32)$,
reversely lavishing $(+14)$ mappings.

The reversely binding $(-21)$ mapping
$$\xymatrix{\centre{-21}
  \cdot 
  & \gimimas  \vent[l]_{f^{P*}} _(1.2)*{\mazaswarrow} \vent[d]^{\Adh_Y}
  \\ \cdot \vent[u]^{\Adh_X^*}
  & \cdot  \vent[l]^-{f^*}}$$
is defined by inclusion of reforms
$$B\circ f^{P*}\subset B\circ\Adh_Y\circ f^*\circ \Adh_X \kbl$$
which is equivalent to the large inverse image of the inverse napping
$$f^{P*}\subset (\Adh_Y\times \Adh_X)\inv (f^*)\tsk$$
This property can also be expressed with union of adherence filters
$$(f^P)\inv (B)\mzl\VF_{f\inv([B])}\kbl$$
every partial set $A\subset X$ which has the image
$f^P(A)=B$ coinciding with the taken
partial set in the target space $B\subset Y$,
has some adherent point in the inverse image $x\in f\inv(B)$
of taken set $B\subset Y$.

\Partnumber,
The thick  $(+43)$ mapping
$$\xymatrix{\centre{+43}
  \cdot \ar[r]^-{f^P} 
  & \cdot \vent[d]^{\Adh_Y}
  \\ \gimimas \vent[u]^{\Adh_X^*} \ar[r]_-{f} _(1.2)*{\mazaswarrow}
   & \cdot }$$
is defined by an inclusion of reforms
$$x\circ \Adh_X^*\circ f^P\circ \Adh_Y\subset x\circ f\kbl$$
which is equivalent for small image of mapping's direct image appointment
$$(\Adh_X\times \Adh_Y)^P(f^P)\subset f\tsk$$
This can be expressed with adherence operator's union
over all sets from some filter
$$\Cup\lbrace [B]:B\in f^{PP}(\VF_x)\rbrace \subset \lbrace f(x)\rbrace\kbl$$
if the taken point $x\in X$ adheres the partial set $A\subset X$,
then the image point $f(x)\in Y$ is a unique adherent point for
the image set $f^P(A)\subset Y$.

For the functional convergence $\Adh_X$ in the source space we get the
equivalent property of wasting mapping $(+13)$. For the functional
convergence $\Adh_Y$ in the target space we get the equivalent property
of creating mapping $(+42)$, therefore for mappings this is equivalent
with reversely respecting $(+31)$, reversely pressing $(+32)$,
reversely lavishing $(+41)$ mappings also.

\Partnumber.
Now rest more special properties, without any equivalent counterparts.
\Partnumber:
The thin $(-12)$ mapping
$$\xymatrix{\centre{-12}
  \gimimas \ar[r]^-{f^P} ^(1.2)*{\mazanwarrow} \vent[d]_{\Adh_X}
  & \cdot 
  \\ \cdot \ar[r]_-{f} 
  & \cdot \vent[u]_{\Adh_Y^*} }$$
is defined by an inclusion  of reforms
$$A\circ \Adh_X\circ f\circ \Adh_Y^* \subset A\circ f^P\kbl$$
which is eqiuvalent to the small inverse image of mapping
$$(\Adh_X\times \Adh_Y)^{*P}(f)\subset f^P \tsk$$
This can be expressed with the union of adherence filters
$$\VF_{f^P([A])}\subset \lbrace f^P(A)\rbrace \kbl$$
the image $f^P(A)\subset Y$ of every partial set $A\subset X$ in
the source space is unique partial set $B\subset Y$ in the target space
which has an adherent point $y\in [B]$ related to some
adherent point of the taken partial set $x\in [A]$, $f(x)=y$.
This means that the related adherence points must also relate
the partial sets.

The reversely thin $(+21)$ mapping
$$\xymatrix{\centre{+21}
  \cdot 
  & \gimimas \vent[d]^{\Adh_y} \vent[l]_{f^{P*}} _(1.2)*{\mazanearrow} 
  \\ \cdot \vent[u]^{\Adh_X^*} 
  & \cdot  \vent[l]^-{f^*} }$$
is defined by inclusion of reforms
$$B\circ \Adh_Y\circ f^*\circ \Adh_X^* \subset B\circ f^{P*}\kbl$$
which is equivalent to the small inverse image of inverse mapping's
 graphic
$$(\Adh_Y\times \Adh_Y)^{*P}(f^*)\subset f^{P*}\tsk$$
This can be expressed with the union of adherence filters
$$\VF_{f\inv ([B])}\subset (f^P)\inv(B)\kbl$$
the partial set $A\subset X$, with the adherent points $x\in [A]$
having related points $y\in Y$, $f(x)=y$ adherent to the taken partial set
$y\in [B]$,  must be related to the taken set in the target space
$f^P(A)=B$.
\Partnumber,
The parting $(-43)$ mapping
$$\xymatrix{\centre{-43}
  \cdot \ar[r]^-{f^P} 
  & \cdot \vent[d]^{\Adh_Y}
  \\ \gimimas \vent[u]^{\Adh_X^*} \ar[r]_-{f} _(1.2)*{\mazanearrow}
  & \cdot  }$$
is defined by inclusion  of reforms
$$x\circ f\subset x\circ\Adh_X\circ f^P\circ \Adh_Y\kbl$$
which is equivalent to a large image of the mapping's
direct image appointment
$$f\subset (\Adh_X\times \Adh_Y)^P(f^P) \tsk$$
The adherence expands set image so much that it include the
 mapping itself $f:X\longrightarrow Y$.
This can be expressed with upper
covering of adherence operator over the adherence filter 
$$\lbrace f(x)\rbrace\subset
\Cup \lbrace [B]:B\in f^{PP}(\VF_x)\rbrace \kbl$$
i. e. the image $f(x)\in Y$ of a taken point $x\in X$ is adherent
to the image $f^P(A)\subset Y$ of any partial set in the source space
$A\subset X$ for which the taken point is adherent $x\in [A]$.

If this adherent point is defined as set's middle point $$x=\int A\kbl$$
then the parting mapping asks that the image $f(x)\in Y$
would be a middle point of an direct image $$f(x)=\int f^P(A)$$
for some set $A\subset X$ for which taken point is middle point.

This property is very week.
For expanding adherences such property has every mapping
$f:X\longrightarrow Y$. We can take $A=\lbrace x\rbrace$, then $x\in [A]$
and $f(x)\in [A]$.

The reversely parting  $(+34)$ mapping
$$\xymatrix{\centre{+34}
  \cdot \vent[d]_{\Adh_X}
  & \cdot \vent[l]_{f^{P*}} 
  \\ \cdot 
  & \gimimas \vent[u]_{\Adh_Y^*} \vent[l]^-{f^*} ^(1.2)*{\mazanwarrow} }$$
is defined by inclusion of reforms
$$y\circ f^*\subset y\circ \Adh_Y^*\circ f^{P*}\circ \Adh_X\kbl$$
which is equivalent to large image of inverse  of mapping's direct image
appointment
$$f^*\subset (\Adh_Y\times \Adh_X)^P (f^{P*})\tsk$$
This can be expressed with upper covering of adherence operator
over the filter
$$f\inv (y)\subset \Cup \lbrace [A]: f^P(A)\in\VF_y\rbrace\kbl $$
i.e. for the point in the target space $y\in Y$ every point of inverce
image  $x\in f\inv(y)$ is adherent some partial set $A\subset X$
which direct image $f^P(A)\subset Y$ is adherent by the taken point
$y\in [f^P(A)]$.
\Partnumber.
We can define weakly cautious mappings of adherence spaces, and
to show that such property is equivalent to respecting mappings of
inside spaces. This provides the fact, that some continuity properties
may be expressed both in adherence and inside spaces. More details
is presented in my recent monograph \nuoroda* Valius: Val>
Valiukevi\v cius G. The continuity property of multivalued functions
(in lithuanian), Veja: Vilnius 2009. *.


%% file: mborno.tex
\Part Mapping's global continuity.
\Partnumber:
The continuity properties in adherence spaces
provide local continuity properties in generated topological
spaces. The same situation also arise  for bornological or measurable
spaces. Therefore it is usefull to look at general continuity properties
for topological spaces. At this moment we content himeself to deal
only with mappings. These continuity properties will be called global.

\Partnumber,
The part in the potential set $\VT\subset X^P$ will be called a
\normalref{=tribe} and it will define a \normalref{=tribe space} $X$.
For any mapping $f: X\rightarrow Y$ with set's direct image appointment
 $f^P:X^P\rightarrow Y^P$ we can define 3 continuity properties:
a \normalref{=carrying mapping} $ f^{PP}(\VT_X)\mzl \VT_Y$,
a \normalref{=hiding mapping} $f^{PP}(\VT_X)\dgl \VT_Y$, and
an \normalref{image =reflecting mapping} $(f^P)\inv(\VT_Y)\mzl \VT_X$.

For topological spaces we can investigate the topology $\VT$ and
cotopology $\VT^c$
together, therefore we get different continuity properties:
  For topology compounded by closed sets
a \normalref{=closed mapping} $ f^{PP}(\VT_X)\mzl \VT_Y$,
a \normalref{closed set =hiding mapping} $f^{PP}(\VT_X)\dgl \VT_Y$, and
an \normalref
{closed sets image =reflecting mapping} $(f^P)\inv(\VT_Y)\mzl \VT_X$,
and for cotopology compounded by open sets
a \normalref{=open mapping} $ f^{PP}(\VT^c_X)\mzl \VT^c_Y$,
a \normalref{open sets =hiding mapping} $f^{PP}(\VT^c_X)\dgl \VT^c_Y$, and
an \normalref{open set image =reflecting mapping} $(f^P)\inv(\VT^c_Y)
\mzl \VT^c_X$.

With set's inverse image appointment  $f\inv: Y^P\rightarrow X^P$
we define another 3 continuity properties:
a \normalref{ =turning mapping} $(f\inv)^P(\VT_Y)\mzl \VT_X$,
a \normalref{=exhausting mapping} $(f\inv)^P(\VT_Y)\dgl \VT_X$,
and an \normalref{=inverse image reflecting mapping}
$(f\inv)\inv (\VT_Y)\mzl \VT_X$.

For topological spaces only the turning mappings are concidered as continuous
ones.

\Partnumber:
We defined the closure as isotonic and expanding adherence. The
closed sets is defined as fixed points for closure operator, i.e.
the closed sets are defined by the property
$$A\subset [A]\subset A\tsk$$
The tribe compounded by all closed sets is called a topology.
It is characterized by the property of freely intersection
$$\Cap \VE\in \VT\uparrow \VE\subset \VT\kbl$$
every tribe, which maintains the arbitrary intersection of its members,
can be get as topology of closed sets for some closure. Such topologies
are used in logic calculations
\nuoroda* Martin: Mar> Martin N. Pollard St.
Closure spaces and logic, Kluwer: Dordrecht -- Boston -- London 1996.*.
Bourbaki for more traditional topologies asks additional property
that wide set is closed.

Otherwise the closure was called a multistep (germ. mehrstufig) topology
by G\"ahler \nuoroda* GahlerGA: Gah> G\"ahler W. Grundstrukturen der Analysis
I, II, Akademie Verlag: Berlin 1977, 1978.*.

The respecting $(-13)$ mapping asks inclusion for direct image
$$f^P([A])\subset [f^P(A)]\uparrow A\subset X\kbl$$
and reversely creating mapping $(-24)$ asks equivalent
(for mappings of closure spaces)
opposite inclusion for inverse image
$$[f\inv(B)]\subset f\inv([B])\uparrow B\subset Y\tsk$$
This can be also expressed by equivalent (for mappings) properties with
direct or inverse images of closure reform. The pressing $(-23)$ mapping
asks an inclusion of direct image
$$(f^P\times f)^P(\Adh_X)\subset Adh_Y\kbl$$
and lavishing $(-14)$ mapping asks the same opposite inclusion of inverse
image
$$\Adh_X\subset (f^P\times f)\inv (\Adh_Y)\tsk$$
These all equivalent continuity properties coincide with usual continuity
of mappings between closure spaces.

Such mapping is turning for topological spaces, i.e. the inverse image of
closed set remains closed
$$[B]\subset B\duoda [f\inv (B)]\subset f\inv([B])\subset f\inv(B)\tsk$$

The wasting $(-31)$ mapping asks opposite inclusion for direct image
$$[f^P(A)]\subset f^P([A])\uparrow A\subset X\tsk$$
It defines the closed mapping between closure spaces:
The direct image of closed set remains closed
 $$[A]\subset A\duoda [f^P(A)]\subset f^P([A])\subset f^P(A)\tsk$$
 Such mappings are carrying for topological spaces.

Bourbaki
 \nuoroda* BourbakiTop: Bou> Bourbaki N. Topologie g\'en\'erale Chapitre 1.
  Structures  topologiques, Hermann: Paris 1961, Nauka: Moskva 1968.*
concidered the property of closed mappings only together with continuity.

The cautious $(-42)$ mappings asks big image of the adherence filter
$$f^{PP}(\VF_x)\dgl \VF_{f(x)}\tsk$$
For covering mapping between closure spaces  it
asks an inclusion of inverse image
$$f\inv([B])\subset [f\inv(B)]\uparrow B\subset Y\tsk$$

The  more useful property of \normalref{weakly =cautious mapping}.
It demands
that image of adherence filter would generate a big hereditary filter
$$[f^{PP}(\VF_x)]\dgl \VF_{f(x)}\tsk$$
Such mapping will be equivalent to respecting mapping of conjugate
interior spaces. It will be open for generated topology.

Let $0\subset X$ is open set of closure in a source space. It is easy
to check that the direct image $f^P(O)\subset Y$ will be open set of
closure in a target space. It is enough to show that any set
$A'\subset Y\setminus f^P(O)$ has no adherent point from image of taken
set $y\in f^P(O)$. Otherwise we should have a point $x\in O$ with $f(x)=y$,
and for weakly cautious mapping this point will be an adhered by the
set $A\subset X$ with smaller direct image $f^P(A)\subset A'$. Therefore
$A\subset X\setminus O$, and such set can't adhere the point from the
taken open set $x\in O$. So we get a contradiction, and there is no point
from the taken open set adhered by the set $A'\subset Y\setminus f^P(O)$.

The hiding $(-32)$ mapping asks
the opposite inclusion of direct image of closure reform
$$(f^P\times f)^P(\Adh_X)\supset \Adh_Y\kbl$$
and reflecting $(+14)$ mapping asks the inclusion of inverse image
 of closure reform
$$(f^P\times f)\inv (\Adh_Y)\subset \Adh_X\tsk$$
Such local properties don't coincide with corresponding global properties
in topological spaces.

\Partnumber,
We define the interior as isotonic and decreasing adherence. The open
sets is defined as fixed points for interior operator, i.e. the open sets are
defined by the property
$$A\supset \langle A\rangle \supset A\tsk$$
The tribe compounded by all open sets we shall call a cotopology. It is
characterized by the property of freely union
$$\Cup\VE\in \VK\uparrow \VE\subset \VK\kbl$$
every tribe, which maintains the arbitrary union of its members,
can be get as cotopology of open sets for some interior. Cotopologies are
dual for topologies, however the continuity properties of mappings
are different.

From continuous mappings of closure sets
with inverse image condition
$$[f\inv(B)]\subset f\inv ([B])\uparrow B\subset Y\kbl$$
we get the dual inverse image condition for mappings of interior spaces
$$f\inv(\langle B\rangle )\subset \langle f\inv (B)\rangle \uparrow B\subset Y\tsk$$
Such mappings are turning between the cotopological spaces
$$B\subset \langle  B\rangle \duoda f\inv (B)\subset f\inv(\langle B\rangle )\subset \langle f\inv(B)\rangle \tsk$$

Nevertheless this condition cannot be get as earlier continuity property
for adherence spaces.
The cautious $(-42)$ mappings between interior spaces is defined by
the same inclusion for covered subsets
$$f\inv(\langle B\rangle )\subset \langle f\inv (B)\rangle \uparrow B\subset f^P(X)\kbl$$
and in the target space we demand that
$$\langle B\rangle \subset Y\setminus f^P(X)\uparrow B\not\subset f^P(X)\tsk$$
Therefore in most cases the cautious mappings must be surjective.

The situation is rescued with the notion of \normalref{weakly =cautious
mapping}
between interior spaces. This property is
equivalent to the property of respecting mapping
for conjugate closure spaces.

The respecting $(-13)$ mapping between interior spaces
is defined by inclusion of direct image
$$f^P(\langle A\rangle )\subset \langle f^P(A)\rangle \uparrow A\subset X\tsk$$
It defines open mapping between topological spaces.
Such mapping is carrying
for cotopological spaces, i.e. the direct image of open set
remains open
$$A\subset \langle A\rangle \duoda f^P(A)\subset f^P(\langle A\rangle )\subset \langle f^P(A)\rangle \tsk$$
Bourbaki in \nuoroda* BourbakiTop: Bou> *
uses such notion together with continuity, ant doesn't remark
any equivalent formulation
of open mappings.

The reversely creating $(-24)$ mapping between interior spaces
is defined by opposite inclusion of inverse image
$$\langle f\inv(B)\rangle \subset f\inv (\langle B\rangle )\uparrow B\subset Y\tsk$$
The pressing $(-23)$ mapping between interior spaces
demands the inclusion of direct image of interior reform
$$(f^P\times f)^P(\Int_X)\subset \Int_Y\tsk$$
The lavishing $(-14)$ mapping between interior spaces
demands the opposite inclusion of inverse image of interior reform
$$\Int_X\subset (f^P\times f)\inv(\Int_Y)\tsk$$
All such mappings define an open mapping of topological spaces.

The wasting $(-31)$ mapping between interior spaces demands the opposite
inclusion of direct image
$$f^P(\langle A\rangle )\supset \langle f^P(A)\rangle
\uparrow A\subset X\tsk$$
This property wasn't remarked for classical topological spaces,
nevertheless it is
meaningful for interesting topologies.

The hiding $(-32)$ mapping between interior spaces demands
the opposite inclusion of direct image of interior reform
$$(f^P\times f)^P(\Int_X)\supset \Int_Y\tsk$$
The reflecting $(+14)$ mapping between interior spaces demands
the the inclusion of inverse image of interior reform
$$(f^P\times f)\inv(\Int_Y)\subset \Int_X\tsk$$

Both last local properties for mappings of cotopological spaces
differ from coresponding global properties.

\Partnumber,
Any tribe $\VK\subset X^P$ can be declared as bornology, and its member
$B\subset X$ as bounded set. The carrying mapping between bornological
spaces $f:X\longrightarrow Y$ is called bounded, i.e.
the image of bounded set must be bounded
$$ f^P(A)\in\VK_Y\uparrow A\in\VK_X\tsk $$
The turning mapping between bornological spaces $f:X\longrightarrow Y$
we shall call perfect, i.e.
the inverse image of bounded set must be bounded
$$f\inv(B)\in\VK_X\uparrow B\in\VK_Y\tsk$$
Such mappings are useful for bornologies compounded of compact parts
in some topological space.

\Partnumber:
It is interesting to consider adherence space as example of "local"
bornology. The reform $\Sid_X:X^P\ventrightarrow X$ will
be called \normalref{=proximity}.
We say that the couple $\langle A,x\rangle \in X^P\times X$ belongs to
proximity  if the set $A$ is
 \normalref{=proximal} to the point $x\in X$.
The direct appointment for each set $A\subset X$
appoints the set of points for
which taken set is proximal, we shall call it a \normalref{=siding set}
with sign $\langle A\rangle \subset X$.
The opposite appointment for every point $x\in X$
appoints a \normalref{=proximal filter} $\VF_x$ of proximal sets

\Partnumber,
In proximity space we define bornology declaring all sets with biger
siding set
$$  A\subset \langle A\rangle\uparrow A\subset X$$
as bounded sets.

In bornology space we define proximity, declaring
the set $A\subset X$ proximal to the point $x\in X$, if we can find
a biger bounded set $A\subset B\subset X$ with taken point $x\in B$.

We shall check continuity properties for mappings between
proximity spaces. Usually we work with proximal cofilters, demanding that
smaller set remainded  bounded and proximal.

\Partnumber,
The respecting $(-13)$ mapping demands inclusion of siding sets
$$f^P(\langle A\rangle)\subset \langle f^P(A)\rangle\tsk$$
Such mappings will be bounded for generated bornology spaces:

Let the set $A\subset X$ is near of every it's point
$$A\subset \langle A\rangle \tsk$$
Then the direct image of taken set $f^P(A)\subset Y$ has the same property
$$f^P(A)\subset f^P(\langle A\rangle )\subset \langle f^P(A)\rangle\tsk$$

And in opposite direction. The bounded mappings of bornological spaces
are respecting for generated proximity spaces:

Let the set $A\subset X$ is near  the point $x\in X$, i. e. we find a bigger
bounded set $A\subset B$ having the taken point $x\in B$. We shall
check that the direct image of taken set $f^P(A)\subset Y$
remains near the point $f(x)\in Y$. It is enough to take bounded set
$f^P(B)\subset Y$, for which
we have $f^P(A)\subset f^P(B)$ and $f(x)\in f^P(B)$.

The reversely respecting $(+31)$ mapping demands
inclusion of proximal filters in source space
$$\VF_{f\inv(y)}\mzl (f^P)\inv(\VF_y)\tsk$$
It is equivalent for creating $(+42)$ mapping, which demands
inclusion of proximal filters in target space
$$f^{PP}(\VF_x)\mzl \VF_{f(x)}\tsk$$
Therefore this property is similar to local boundedness.

There are additional two equivalent properties defined with direct or
inverse images of the proximity reform.

The reversely pressing $(+32)$ mapping between proximity spaces
demands a small direct image
$$ (f^P\times f)^P(\Sid_X^*)\subset \Sid_Y\kbl$$
and reversely lavishing $(+41)$ mapping between proximity spaces
demands a big inverse image
$$ \Sid_X^*\subset (f^P\times f()\inv (\Sid_Y^*)\tsk$$

\Partnumber,
Rest 4 more special continuity properties.

The  cautious $(-42)$ mapping between proximity spaces
demands big direct image of proximity filter
$$f^{PP}(\VF_x)\dgl \VF_{f(x)}\tsk$$
We can show that such mappings are perfect for generated bornologies:
The inverse image $f\inv(B')\subset X$ of bounded set $B'\subset Y$ remains
bounded.

Let the set $B'\subset Y$ has a bigger siding set
$B'\subset \langle B'\rangle$. We shall check that the same is also an
 inverse image
$$f\inv(B')\subset \langle f\inv(B')\rangle \tsk$$
The smaller set $B''\subset B'$ remains siding, therefore a smaller set
remains bounded
$$ B''\subset B'\subset \langle B'\rangle \subset \langle B''\rangle\tsk$$
We can take a smaller bounded set $B'':=f^P(f\inv(B'))\subset B'$ covered
by mapping,
therefore for cautious mapping we get
$$f\inv(\langle B''\rangle)\subset \langle f\inv(B'')\rangle \kbl $$
$$f\inv(B')=f\inv(B'')\subset f\inv(\langle B''\rangle \subset
\langle f\inv(B'')\rangle =\langle f\inv(B')\rangle\tsk$$

And in opposite direction. For perfect covering mapping we get reversely
cautious mapping
$(+24)$ of generated proximity spaces, i. e. for the set $A'\subset Y$
the inverse image of siding set $f\inv(\langle A'\rangle )\subset X$
must be included in union of siding sets
$$f\inv (\langle A'\rangle)\subset
\Cup \lbrace \langle A\rangle: f^P(A)=A'\rbrace \tsk$$
Every siding point $y\in \langle A'\rangle$ is included in bigger bounded set
$y\in B'$, $A'\subset B'$. We can take the inverse image of this bounded set
$f\inv(y)\subset f\inv(B')$ and it will be bigger for any set $A\subset X$
with $f^P(A)=A'$. For covering mapping we can take $A:=f\inv(A')$, therefore
$$f\inv(\langle A'\rangle)\subset \langle f\inv (A')\rangle\tsk$$

The reflecting mappings $(+14)$ demmands small inverse image of siding
reform
$$(f^P\times f)\inv(\Sid_Y)\subset \Sid_X\tsk$$
We shall check that it reflects bounded sets for generated bornologies.

Let we have the set $A\subset X$ with bounded direct image
$$f^P(A)\subset \langle f^P(A)\rangle \tsk$$
We shall show that taken set is also bounded
$$A\subset \langle A\rangle\tsk$$
Every point $x\in A$ has image $f(x)\in Y$ proximal to the direct image of
taken set $f^P(A)\subset Y$, therefore for reflecting mapping it will
be proximal to the taken set $x\in \langle A\rangle $.

And in oposite direction. The mapping $f: X\rightarrow Y$ which
reflects bounded sets will be also reflecting for generated proximities.

Let $f(x)\in Y$ is siding to the set $f^P(A)\subset Y$, i. e. we
can find a bigger bounded set $f^P(A)\subset B'$ with $f(x)\in B'$.
We can take $B:=f\inv(B')$, it's image remains bounded
$f^P(A)\subset f^P(B) \subset B'$, therefore for reflecting mapping it itself
is bounded. We have find a bigger bounded set $A\subset B$ with $x\in B$.

The hiding $(+23)$ mapping between siding spaces
demands large direct image of siding reform
$$\Sid_Y\subset (f^P\times f)^P(\Sid_X)\tsk$$
We shall check that all mappings hiding the bounded sets of bornology spaces
are also hiding between siding spaces.

Let we have the point $y\in Y$ and set $A'\subset Y$ siding in the
bornology space $Y$, i. e. we can find a bigger bounded set $A'\subset B'$
with taken point $y\in B'$. For hiding mapping we can find a bounded set
$B\subset X$ with direct image $f^P(B)=B'$. Therefore we also can find
a covering point $x\in X$ with $f(x)=y$. Taking the set $A:=f\inv(A')\cap B$
we get needed siding set $f^P(A)=A'$.

In general we can't check that hiding  mapping between siding spaces
is hiding the bounded sets.

At first we notice that every bounded point $y\in Y$
must be covered. Such point has proximal set $A'\subset Y$, therefore for
hiding mapping
will be the point $x\in X$ with $f(x)=y$ and the set $A\subset X$ with
$f^P(A)=B$.

Let the set $B'\subset Y$ is bounded for the siding space, i. e.
it has a bigger siding set
$B'\subset \langle B'\rangle$. Every point $y\in B'$ will have $x\in X$
with $f(x)=y$
and $B_x\subset X$ with $f^P(B_x)=B'$, therefore these points compounds the
set $A\subset X$ with direct image $f^P(A)=B'$. We want to check that the
new set has bigger siding set
$$A\subset \langle A\rangle $$
But it can't be calculated without some asumption of compactness
for the union of proximal sets
$$\Cup \lbrace B_x: x\in X\rbrace \tsk$$

The reversely wasting $(+13)$ mapping between siding spaces
demands small inverse image of siding filter
$$(f^P)\inv (\VF_y)\mzl \VF_{f\inv(y)}\tsk$$
I don't think that these property may have global counterpart for bornology
spaces.

At this moment I don't see any application of such continuity properties
for bornological spaces. May be they become interesting for
various final or initial constructions of siding spaces.
